\documentclass[10pt]{article}

\usepackage{graphics}
\usepackage{graphicx}

\usepackage[a4paper, left=35mm,right=35mm,top=34mm,bottom=34mm]{geometry}
\usepackage[utf8]{inputenc}
\usepackage[T1]{fontenc}
\usepackage[english]{babel}

\usepackage{enumerate}
\usepackage{graphicx}
\usepackage{hyperref}
\hypersetup{
    colorlinks=true,
    linkcolor=blue,
    filecolor=magenta,      
    urlcolor=cyan,
}
\usepackage{listings}
\usepackage{color}

\definecolor{dkgreen}{rgb}{0,0.6,0}
\definecolor{gray}{rgb}{0.5,0.5,0.5}
\definecolor{mauve}{rgb}{0.58,0,0.82}
\lstdefinelanguage{MRGC++}{%
  language=C++,
  morekeywords={T, U, MPI_Irecv, MPI_Isend, MPI_Allreduce, MPI_Waitall, Compute, Map, abs, max, Swap, MPI_Recv_init, MPI_Send_init, MPI_Startall, Copy, Init, InitRecv, InitSend, InitAllReduce, Send, Recv, AllReduce, Finalize, InitSnapshot, Snapshot, SwitchAsync, SnapReduce, MPI_Test, MPI_Start}
}
\usepackage{mathtools,amsthm,amssymb,amsfonts}
\usepackage{algorithm}
\usepackage{algorithmic}
\makeatother
\theoremstyle{plain}
\newtheorem{theorem}{Theorem}

\newtheorem{lemma}{Lemma}
\theoremstyle{definition}

\theoremstyle{definition}

\theoremstyle{remark}
\newtheorem*{remark}{Remark}

\usepackage{caption} 
\usepackage{fancyhdr}

\author{
  {\normalsize Guillaume Gbikpi-Benissan}\thanks{Universit\'e Paris-Saclay, CentraleSup\'elec, France.}
  \and
  {\normalsize Fr\'ed\'eric Magoul\`es}\thanks{Universit\'e Paris-Saclay, CentraleSup\'elec, France \textbf{/} Faculty of Engineering and Information Technology, University of P\'ecs, Hungary
    (correspondence, frederic.magoules@hotmail.com).}
}
\title{Asynchronous multiplicative coarse-space correction}
\date{}

\begin{document}
\maketitle
\thispagestyle{fancy}

\begin{abstract}
\noindent This paper introduces the multiplicative variant of the recently proposed asynchronous additive coarse-space correction method.
Definition of an asynchronous extension of multiplicative correction is not straightforward, however, our analysis allows for usual asynchronous programming approaches. General asynchronous iterative models are explicitly devised both for shared or replicated coarse problems and for centralized or distributed ones. Convergence conditions are derived and shown to be satisfied for M-matrices, as also done for the additive case.
Implementation aspects are discussed, which reveal the need for non-blocking synchronization for building the successive right-hand-side vectors of the coarse problem. Optionally, a parameter allows for applying each coarse solution a maximum number of times, which has an impact on the algorithm efficiency. Numerical results on a high-speed homogeneous cluster { confirm} the practical efficiency of the asynchronous two-level method over its synchronous counterpart, even when it is not the case for the underlying one-level methods.
\end{abstract}


\section{Introduction}

We are considering parallel iterative solution of sparse linear problems
\begin{equation}
\label{eq:axeb}
A x = b, \qquad A \in \mathbb{R}^{n \times n},
\end{equation}
arising, e.g., from finite-elements or finite-differences discretization of partial differential equations. One of the most popular techniques to achieve faster and scaling parallel iterative methods consists of correcting the computed approximate solution $x^{k}$ at iteration $k$ such as
\[
\bar x^{k} := x^{k} + P \widetilde y^{k},
\]
where $\widetilde y^{k}$ is the solution of a coarser problem,
\[
\widetilde A \widetilde y^{k} = R \left(b - A x^{k}\right).
\]
$P$ and $R$ are, respectively, prolongation and restriction operators for mapping between the fine and the coarse spaces. Iterative methods under consideration here are of the form
\begin{equation}
\label{eq:si}
x^{k+1} = x^{k} + M \left(b - A x^{k}\right),
\end{equation}
where $M$ is an approximation to the inverse of $A$. Such a technique can be recursively applied to the coarser problem as well, leading to a multilevel method. We refer the reader to, e.g., \cite{Hackbusch1985} for a comprehensive introductory review about principles of multilevel methods and to \cite{TosWid2005} for their developments within domain decomposition frameworks.

Our specific interest here is the extension of such multilevel methods to asynchronous iterative models \cite{Baudet1978}. Asynchronous iterations were first experimented in \cite{Rosenfeld1969} as a generalization of free steering techniques (see, e.g., \cite{Schechter1959}), leading to the development of an extensive general convergence theory (see, e.g., \cite{ChazMir1969, Baudet1978, Bertsekas1983, FromSzyld1994}). Quite early, asynchronous multilevel methods were investigated as in \cite{HartCormick1989} where slight modifications were operated on the fast adaptive composite grid method to remove sequential dependence between the refinement levels. Latest related works include, e.g., the asynchronous Parareal time integration scheme \cite{MagGBen2018} which is a two-level time-parallel method. In \cite{HawkesEtAl2019}, multilevel coarse solutions are asynchronously computed and propagated, while synchronization still occurs at the beginning of each cycle. Finally, asynchronous additive coarse-space correction was proposed in \cite{WPouChow2019, GlusaEtAl2020} for the restricted additive Schwarz method.
An alternative approach is proposed in this paper, where asynchronous multiplicative coarse-space correction is modeled and analyzed within a general iterative framework including additive Schwarz-type methods.

The paper is organized as follows. Section \ref{sec:bg} recalls notions of Schwarz-type methods, asynchronous iterations and multiplicative coarse-space correction. Section \ref{sec:theo} presents our theoretical results about convergence conditions of two-level asynchronous iterative methods, depending on two different parallel schemes. Section \ref{sec:impl} discusses implementation aspects along with detailed examples. Section \ref{sec:exp} presents our experimental results about optimal parameters and performance scaling, { using a three-dimensional Poisson's problem}. Our conclusions follow in Section \ref{sec:concl}.

\section{Computational background}
\label{sec:bg}

\subsection{Asynchronous iterative methods}

Let $\Omega$ { denote} the set of $n$ unknowns from the problem \eqref{eq:axeb} and assume {
a possibly overlapping} distribution into $p$ subsets $\Omega^{(1)}$ to $\Omega^{(p)}$ with $p \le n$. Induce, then, $p$ square sub-matrices $A^{(1)} \in \mathbb{R}^{n^{(1)} \times n^{(1)}}$ to $A^{(p)} \in \mathbb{R}^{n^{(p)} \times n^{(p)}}$ such that
{
\[
A^{(s)} = R^{(s)} A {R^{(s)}}^{\mathsf T}, \qquad b^{(s)} = R^{(s)} b, \qquad s \in \{1, \ldots, p\},
\]
where $R^{(s)}$ is the unknowns selection mapping from $\Omega$ to $\Omega^{(s)}$. Assume, now, diagonal weighting matrices $W^{(1)}$ to $W^{(p)}$ such that
\begin{equation}
\label{eq:weighting_constraint}
\sum_{s=1}^{p} {R^{(s)}}^{\mathsf T} W^{(s)} R^{(s)} = I,
\end{equation}
and consider, then, the classical iterative scheme \eqref{eq:si} with
\[
M = \sum_{s=1}^{p} {R^{(s)}}^{\mathsf T} W^{(s)} M^{(s)} R^{(s)}.
\]
}
{
For the unknowns not in the overlap part of $\Omega^{(s)}$, the corresponding diagonal entries in $W^{(s)}$ are obviously $1$. Differences in weighting strategies therefore occur in
\[
\Omega^{(s)} \cap \bigcup_{\substack{r=1\\ r \ne s}}^{p} \Omega^{(r)}.
\]
}
Note, first, that if $M^{(s)} = {A^{(s)}}^{-1}$, then $M$ is a Schwarz-type preconditioner. If, then, $W^{(s)} = I^{(s)}$ (identity on $\Omega^{(s)}$), then the constraint \eqref{eq:weighting_constraint} induces { no overlap (or minimum overlap at continuous level)} and so, $M$ results in a block-Jacobi preconditioner. However, simply removing $W^{(s)}$ from the expression of $M$ would give the additive Schwarz (AS) preconditioner (with arbitrary overlap). If the entries of $W^{(s)}$ rather consist of $0$ and $1$, then the constraint \eqref{eq:weighting_constraint} induces a restricted AS preconditioner. Allowing, finally, for nonzero weights less than $1$ (typically, $1/m$ with $m$ being the number of subdomains sharing the unknown), we obtain a weighted restricted AS preconditioner. The reader is referred to, e.g, the original papers \cite{FrommerSchwandt1997, CaiSarkis1999} and the review paper \cite{Gander2008}, for detailed insights into Schwarz-type methods.

Now, from \eqref{eq:weighting_constraint}, we have, for any vector $x \in \mathbb{R}^{n}$,
\begin{equation}
\label{eq:vec_decomp}
x = \sum_{s=1}^{p} {R^{(s)}}^{\mathsf T} W^{(s)} {x^{(s)}}, \qquad x^{(s)} := R^{(s)} x,
\end{equation}
hence, the iterative scheme \eqref{eq:si} develops into
\[
\begin{aligned}
\sum_{s=1}^{p} {R^{(s)}}^{\mathsf T} W^{(s)} {x^{(s)}}^{k+1} & = \sum_{s=1}^{p} {R^{(s)}}^{\mathsf T} W^{(s)}\left[{x^{(s)}}^{k} + M^{(s)} R^{(s)} \left(b - A x^{k}\right)\right],
\end{aligned}
\]
which allows for the parallel scheme
\[
\begin{aligned}
{x^{(s)}}^{k+1} & = {x^{(s)}}^{k} + M^{(s)} R^{(s)} \left(b - A x^{k}\right)\\
& = {x^{(s)}}^{k} + M^{(s)} \left(b^{(s)} - \sum_{r=1}^{p} R^{(s)} A {R^{(r)}}^{\mathsf T} W^{(r)} {x^{(r)}}^{k}\right).
\end{aligned}
\]
{
A corresponding computational model is described in Algorithm \ref{algo:si}.
\begin{algorithm}[htbp]
\caption{Iterative solver}
\label{algo:si}
{
\begin{algorithmic}[1]
\STATE \textbf{initialize} local solution $x^{(s)}$ and interface data $x^{(r)}$, $r \ne s$
\WHILE{no convergence}
	\STATE{$x^{(s)}$ := ${x^{(s)}} + M^{(s)} \left(b^{(s)} - \displaystyle\sum_{r=1}^{p} R^{(s)} A {R^{(r)}}^{\mathsf T} W^{(r)} {x^{(r)}}\right)$}
	\STATE \textbf{send} interface data from $x^{(s)}$ \textbf{to} processes $r \ne s$
	\STATE \textbf{receive} interface data $x^{(r)}$ \textbf{from} processes $r \ne s$
\ENDWHILE
\end{algorithmic}
}
\end{algorithm}
}

The generalization into an asynchronous scheme consists of allowing each sub-vector { ${x^{(s)}}$ to be updated without waiting for} synchronization with the other sub-vectors. This is modeled by considering, at each iteration $k$, an arbitrary set $S_{k} \subseteq \{1, \ldots, p\}$ indicating sub-vectors which are updated at that moment. The absence of synchronization further implies to potentially use an outdated component ${x^{(r)}}^{\delta_{s,r}(k)}$ from an arbitrary previous iteration $\delta_{s,r}(k) \le k$. In practice, this models delays in access of a process $s$ to data computed by a process $r$. We shall however naturally require that any $s \in \{1, \ldots, p\}$ belongs to $S_{k}$ an infinite number of times, so that computation of ${x^{(s)}}^{k+1}$ never stops as $k$ grows. Similarly, any $\delta_{s,r}(k)$ is an unbounded increasing function, so that ${x^{(r)}}^{\delta_{s,r}(k)}$ is continuously newer as $k$ grows. Asynchronous iterative methods for linear algebraic problems are thus given by
\begin{equation}
\label{eq:ai}
{x^{(s)}}^{k+1} = \left \{
\begin{array}{ll}
f^{(s)}\left({x^{(1)}}^{\delta_{s,1}(k)}, \ldots, {x^{(p)}}^{\delta_{s,p}(k)}\right) & \forall s \in S_{k},\\
{x^{(s)}}^{k} & \forall s \notin S_{k}
\end{array}
\right.
\end{equation}
with the requirements
\[
\operatorname{card}\left\{k \in \mathbb{N} \ | \ s \in S_{k}\right\} = +\infty,
\qquad
\lim_{k \to +\infty} \delta_{s,r}(k) = +\infty,
\]
where the mapping $f^{(s)}$ is defined as
\[
f^{(s)}\left({x^{(1)}}, \ldots, {x^{(p)}}\right) := {x^{(s)}} + M^{(s)} \left(b^{(s)} - \sum_{r=1}^{p} R^{(s)} A {R^{(r)}}^{\mathsf T} W^{(r)} {x^{(r)}}\right).
\]

\subsection{Coarse-space correction}

Consider, again, the iterative scheme \eqref{eq:si}. Given an iterate $x^{k}$, one would like to correct it with a vector $y^{k}$ such that
\[
A\left(x^{k}+y^{k}\right) = b,
\]
which results in a problem
\[
A y^{k} = b - A x^{k},
\]
however, of the same order as \eqref{eq:axeb}. The idea of multilevel methods therefore consist of instead solving a coarse, cheaper, version of the correction problem,
\begin{equation}
\label{eq:coarse_pb}
\widetilde A \widetilde y^{k} = \widetilde R \left(b - A x^{k}\right),
\end{equation}
where $\widetilde R$ is a restriction mapping from the fine to the coarse space. The correction is then applied as
\[
x^{k} + \widetilde P \widetilde y^{k},
\]
where $\widetilde P$ is a prolongation mapping from the coarse to the fine space. It yields the two-level iterative scheme
\[
\begin{aligned}
x^{k+1} & = x^{k} + \widetilde P \widetilde y^{k} + M \left(b - A \left(x^{k} + \widetilde P \widetilde y^{k}\right)\right)\\
& = x^{k} + \widetilde P \widetilde A^{-1} \widetilde R \left(b - A x^{k}\right) + M \left(b - A \left(x^{k} + \widetilde P \widetilde A^{-1} \widetilde R \left(b - A x^{k}\right)\right)\right)\\
& = x^{k} + \left(M + \widetilde P \widetilde A^{-1} \widetilde R - M A \widetilde P \widetilde A^{-1} \widetilde R\right)\left(b - A x^{k}\right),
\end{aligned}
\]
which exhibits a preconditioner
\[
M + N - M A N, \qquad N := \widetilde P \widetilde A^{-1} \widetilde R,
\]
and an iteration matrix
\[
I - \left(M + N - M A N\right) A = \left(I - M A\right)\left(I - N A\right).
\]
The alternating form of such two-level methods can be written as
\[
\left \{
\begin{array}{lcl}
x^{k+1/2} & = & x^{k} + N \left(b - A x^{k}\right),\\
x^{k+1} & = & x^{k+1/2} + M \left(b - A x^{k+1/2}\right),
\end{array}
\right.
\]
which develops into the parallel scheme
\begin{equation}
\label{eq:csi}
\left \{
\begin{array}{lcl}
{x^{(s)}}^{k+1/2} & = & {x^{(s)}}^{k} + R^{(s)} \widetilde P \widetilde A^{-1} \left(\widetilde b - \displaystyle\sum_{r=1}^{p} \widetilde R A {R^{(r)}}^{\mathsf T} W^{(r)} {x^{(r)}}^{k}\right),\\
{x^{(s)}}^{k+1} & = & f^{(s)}\left({x^{(1)}}^{k+1/2}, \ldots, {x^{(p)}}^{k+1/2}\right)
\end{array}
\right.
\end{equation}
with $s \in \{1, \ldots, p\}$.
On implementation aspects, this corresponds to a scheme where the whole coarse problem is accessed by each of the computing processes, either by memory sharing or by replication.
{
A corresponding replication-based computational model is described in Algorithm \ref{algo:csi}.
\begin{algorithm}[htbp]
\caption{Two-level iterative solver with replicated coarse problem}
\label{algo:csi}
{
\begin{algorithmic}[1]
\STATE \textbf{initialize} local solution $x^{(s)}$ and interface data $x^{(r)}$, $r \ne s$
\WHILE{no convergence}
	\STATE $\widetilde \tau^{(s)}$ := $\widetilde R {R^{(s)}}^{\mathsf T} W^{(s)} \left(b^{(s)} - \displaystyle\sum_{r=1}^{p} R^{(s)} A {R^{(r)}}^{\mathsf T} W^{(r)} {x^{(r)}}\right)$
	\STATE \textbf{send} coarse vector component $\widetilde \tau^{(s)}$ \textbf{to} processes $r \ne s$
	\STATE \textbf{receive} coarse vector components $\widetilde \tau^{(r)}$ \textbf{from} processes $r \ne s$
	\STATE{$\widetilde x$ := $\widetilde A^{-1} \displaystyle\sum_{r=1}^{p} \widetilde \tau^{(r)}$}
	\STATE{$x^{(s)}$ := $x^{(s)} + R^{(s)}\widetilde P \widetilde x$}
	\STATE{$x^{(r)}$ := $x^{(r)} + R^{(r)}\widetilde P \widetilde x$ \COMMENT{correction of interface data}}
	\STATE{$x^{(s)}$ := ${x^{(s)}} + M^{(s)} \left(b^{(s)} - \displaystyle\sum_{r=1}^{p} R^{(s)} A {R^{(r)}}^{\mathsf T} W^{(r)} {x^{(r)}}\right)$}
	\STATE \textbf{send} interface data from $x^{(s)}$ \textbf{to} processes $r \ne s$
	\STATE \textbf{receive} interface data $x^{(r)}$ \textbf{from} processes $r \ne s$
\ENDWHILE
\end{algorithmic}
}
\end{algorithm}
}
A centralized or distributed version would rather be given by
\begin{equation}
\label{eq:csi2}
\left \{
\begin{array}{lcl}
\widetilde y^{k} & := & \widetilde A^{-1} \left(\widetilde b - \displaystyle\sum_{r=1}^{p} \widetilde R A {R^{(r)}}^{\mathsf T} W^{(r)} {x^{(r)}}^{k}\right),\\
{x^{(s)}}^{k+1/2} & = & {x^{(s)}}^{k} + R^{(s)} \widetilde P \widetilde y^{k},\\
{x^{(s)}}^{k+1} & = & f^{(s)}\left({x^{(1)}}^{k+1/2}, \ldots, {x^{(p)}}^{k+1/2}\right),
\end{array}
\right.
\end{equation}
where a process or a group of processes is dedicated to the evaluation of the coarse solution.
{
A corresponding centralization-based computational model is described in Algorithm \ref{algo:csi2}.
\begin{algorithm}[htbp]
\caption{Two-level iterative solver with centralized coarse problem}
\label{algo:csi2}
{
\begin{algorithmic}[1]
\STATE \textbf{initialize} local solution $x^{(s)}$ and interface data $x^{(r)}$, $r \ne s$
\WHILE{no convergence}
	\STATE $\widetilde \tau^{(s)}$ := $\widetilde R {R^{(s)}}^{\mathsf T} W^{(s)} \left(b^{(s)} - \displaystyle\sum_{r=1}^{p} R^{(s)} A {R^{(r)}}^{\mathsf T} W^{(r)} {x^{(r)}}\right)$
	\IF{$s$ = $s_{0}$}
		\STATE \textbf{receive} coarse vector components $\widetilde \tau^{(r)}$ \textbf{from} processes $r \ne s$
		\STATE{$\widetilde x$ := $\widetilde A^{-1} \displaystyle\sum_{r=1}^{p} \widetilde \tau^{(r)}$}
		\STATE \textbf{send} coarse solution $\widetilde x$ \textbf{to} processes $r \ne s$
	\ELSE
		\STATE \textbf{send} coarse vector component $\widetilde \tau^{(s)}$ \textbf{to} process $s_{0}$
		\STATE \textbf{receive} coarse solution $\widetilde x$ \textbf{from} process $s_{0}$
	\ENDIF
	\STATE{$x^{(s)}$ := $x^{(s)} + R^{(s)}\widetilde P \widetilde x$}
	\STATE{$x^{(r)}$ := $x^{(r)} + R^{(r)}\widetilde P \widetilde x$ \COMMENT{correction of interface data}}
	\STATE{$x^{(s)}$ := ${x^{(s)}} + M^{(s)} \left(b^{(s)} - \displaystyle\sum_{r=1}^{p} R^{(s)} A {R^{(r)}}^{\mathsf T} W^{(r)} {x^{(r)}}\right)$}
	\STATE \textbf{send} interface data from $x^{(s)}$ \textbf{to} processes $r \ne s$
	\STATE \textbf{receive} interface data $x^{(r)}$ \textbf{from} processes $r \ne s$
\ENDWHILE
\end{algorithmic}
}
\end{algorithm}
}
While only implementation features make us distinguish \eqref{eq:csi} from \eqref{eq:csi2}, we shall see that their asynchronous versions do result in two mathematically different iterative models.
{
\begin{remark}
It may happen that one applies $f^{(s)}$ before the coarse-space correction. In such case, since these two operations alternate, the two-level schemes \eqref{eq:csi} and \eqref{eq:csi2} remain applicable by considering ${x^{(s)}}^{1/2}$ (the first application of $f^{(s)}$) as the initial guess instead of ${x^{(s)}}^{0}$.
\end{remark}
}

\section{Theoretical results}
\label{sec:theo}

\subsection{Preliminaries}

By $\rho(A)$ and $|A|$, we denote, respectively, the spectral radius and the absolute value (entry-wise) of a matrix $A$. By $\max(x, y)$, we denote an entry-wise maximum operator which selects the maximum of each corresponding entries of vectors $x$ and $y$. By $\|A\|_{\infty}^{w}$, we denote a weighted maximum norm of a matrix $A = (a_{i,j})$, defined as
\[
\|A\|_{\infty}^{w} := \max_{i} \frac{1}{w_{i}} \sum_{j} |a_{i,j}| w_{j}, \qquad w > 0,
\]
where $w$ is a given positive vector. By $\mathsf M$-matrix, we designate a nonsingular matrix with no positive off-diagonal entry which has a nonnegative inverse.

It is well established that one-level iterative methods \eqref{eq:si} are convergent if and only if $\rho(I-MA) < 1$ (see, e.g, Theorem 2.1 in \cite{BahiEtAl2007}). Similarly, one-level asynchronous iterative methods \eqref{eq:ai} were shown in \cite{ChazMir1969} to be convergent if and only if $\rho(|I-MA|) < 1$. Analogously, two-level iterative methods \eqref{eq:csi} are convergent if and only if $\rho((I-MA)(I-NA)) < 1$, where $N$ is the preconditioner induced by the coarse solution. However, asynchronous convergence condition of such methods cannot be trivially inferred as $\rho(|(I-MA)(I-NA)|) < 1$ since the iteration matrix $(I-MA)(I-NA)$ only results from synchronizing $x^{k+1/2}$ and $x^{k+1}$. In the following section, we shall however derive a quite close sufficient condition for the convergence of two-level asynchronous iterations with multiplicative coarse-space correction.

\subsection{Shared/replicated coarse problem}

To enable asynchronous formulation, let us first rewrite the alternating scheme \eqref{eq:csi} as
\[
\left \{
\begin{array}{lcl}
{\chi^{(s)}}^{k} & := & {x^{(s)}}^{k} + R^{(s)} \widetilde P \widetilde A^{-1} \left(\widetilde b - \displaystyle\sum_{r=1}^{p} \widetilde R A {R^{(r)}}^{\mathsf T} W^{(r)} {x^{(r)}}^{k}\right),\\
{x^{(s)}}^{k+1} & = & f^{(s)}\left({\chi^{(1)}}^{k}, \ldots, {\chi^{(p)}}^{k}\right),
\end{array}
\right.
\]
which now can take the asynchronous form
\begin{equation}
\label{eq:cai}
\left \{
\begin{array}{lcl}
{\chi^{(s)}}^{k} & := & {x^{(s)}}^{\delta_{s,s}(k)}\\
&& +\ R^{(s)} \widetilde P \widetilde A^{-1} \left(\widetilde b - \displaystyle\sum_{r=1}^{p} \widetilde R A {R^{(r)}}^{\mathsf T} W^{(r)} {x^{(r)}}^{\delta_{s,r}(k)}\right),\\
{x^{(s)}}^{k+1} & = & \left\{
\begin{array}{ll}
f^{(s)}\left({\chi^{(1)}}^{\delta_{s,1}(k)}, \ldots, {\chi^{(p)}}^{\delta_{s,p}(k)}\right) & \forall s \in S_{k},\\
{x^{(s)}}^{k} & \forall s \notin S_{k},
\end{array}
\right.
\end{array}
\right.
\end{equation}
where each process $s \in \{1, \ldots, p\}$ handles a proper coarse right-hand-side vector
\[
\widetilde b - \displaystyle\sum_{r=1}^{p} \widetilde R A {R^{(r)}}^{\mathsf T} W^{(r)} {x^{(r)}}^{\delta_{s,r}(k)},
\]
leading to multiple, potentially different, coarse solutions. It clearly appears that for each ${\chi^{(r)}}^{\delta_{s,r}(k)}$ with $r \in \{1, \ldots, p\}$, a specific set of sub-vectors ${x^{(q)}}^{\delta_{r,q}(\delta_{s,r}(k))}$ with $q \in \{1, \ldots, p\}$ is involved in the computation of ${x^{(s)}}^{k+1}$ with $s \in S_{k}$. The resulting iteration mapping, say, $F^{(s)}$ thus operates on $p$ versions, $\delta_{1,q}(\delta_{s,1}(k))$ to $\delta_{p,q}(\delta_{s,p}(k))$, of each sub-vector $x^{(q)}$, which does not fit the usual asynchronous scheme \eqref{eq:ai}, where only one version of each sub-vector is involved. While an iteration matrix $(I-MA)(I-NA)$ with $N := \widetilde P \widetilde A^{-1} \widetilde R$ cannot be derived here, a related condition still ensures convergence as it follows.
\begin{theorem}
\label{theo:cai_conv}
A two-level asynchronous iterative method \eqref{eq:cai} is convergent if
\begin{equation}
\label{eq:cai_conv_cond}
\rho\left(\sum_{s=1}^{p}\left|(I - M A) {R^{(s)}}^{\mathsf T} W^{(s)} R^{(s)} \left(I - N A\right)\right|\right) < 1, \qquad N := \widetilde P \widetilde A^{-1} \widetilde R.
\end{equation}
\end{theorem}
\begin{proof}
The model \eqref{eq:cai} develops, for $s \in S_{k}$, as
\[
\begin{aligned}
{x^{(s)}}^{k+1} & = {\chi^{(s)}}^{\delta_{s,s}(k)} + M^{(s)} \left(b^{(s)} - \sum_{r=1}^{p} R^{(s)} A {R^{(r)}}^{\mathsf T} W^{(r)} {\chi^{(r)}}^{\delta_{s,r}(k)}\right)\\
& = {x^{(s)}}^{\delta_{s,s}(\delta_{s,s}(k))} + R^{(s)} \widetilde P \widetilde A^{-1} \left(\widetilde b - \displaystyle\sum_{q=1}^{p} \widetilde R A {R^{(q)}}^{\mathsf T} W^{(q)} {x^{(q)}}^{\delta_{s,q}(\delta_{s,s}(k))}\right)\\
& \quad + M^{(s)} \left(b^{(s)} - \sum_{r=1}^{p} R^{(s)} A {R^{(r)}}^{\mathsf T} W^{(r)}\left( {x^{(r)}}^{\delta_{r,r}(\delta_{s,r}(k))} \right.\right.\\
& \left.\left. \qquad\qquad\quad + R^{(r)} \widetilde P \widetilde A^{-1} \left(\widetilde b - \displaystyle\sum_{q=1}^{p} \widetilde R A {R^{(q)}}^{\mathsf T} W^{(q)} {x^{(q)}}^{\delta_{r,q}(\delta_{s,r}(k))}\right)\right)\right).
\end{aligned}
\]
Taking, then, $p$ vectors $x_{1}$ to $x_{p}$ in $\mathbb{R}^{n}$, the corresponding iteration mapping $F^{(s)}$ is given by
{
\[
\begin{aligned}
F^{(s)}\left(x_{1}, \ldots, x_{p}\right) & := {x_{s}}^{(s)} + R^{(s)} \widetilde P \widetilde A^{-1} \left(\widetilde b - \displaystyle\sum_{q=1}^{p} \widetilde R A {R^{(q)}}^{\mathsf T} W^{(q)} {x_{s}}^{(q)}\right)\\
& \quad + M^{(s)} \left(b^{(s)} - \sum_{r=1}^{p} R^{(s)} A {R^{(r)}}^{\mathsf T} W^{(r)}\left( {x_{r}}^{(r)} \right.\right.\\
& \left.\left. \qquad\qquad\quad + R^{(r)} \widetilde P \widetilde A^{-1} \left(\widetilde b - \displaystyle\sum_{q=1}^{p} \widetilde R A {R^{(q)}}^{\mathsf T} W^{(q)} {x_{r}}^{(q)}\right)\right)\right).
\end{aligned}
\]
}
From \eqref{eq:vec_decomp}, it yields
\[
\begin{aligned}
F^{(s)}\left(x_{1}, \ldots, x_{p}\right) & = R^{(s)} {x_{s}} + R^{(s)} \widetilde P \widetilde A^{-1} \widetilde R \left(b - A {x_{s}}\right)\\
& \quad + M^{(s)} R^{(s)} \left(b - A \sum_{r=1}^{p} {R^{(r)}}^{\mathsf T} W^{(r)}\left(R^{(r)} {x_{r}} \right.\right.\\
& \left.\left. \qquad\qquad\qquad\quad +\ R^{(r)} \widetilde P \widetilde A^{-1} \widetilde R \left(b - A {x_{r}}\right)\right)\right).
\end{aligned}
\]
The induced global iteration mapping follows as
\[
\begin{aligned}
F\left(x_{1}, \ldots, x_{p}\right) & = \sum_{s=1}^{p} {R^{(s)}}^{\mathsf T} W^{(s)} F^{(s)}\left(x_{1}, \ldots, x_{p}\right)\\
& = \sum_{s=1}^{p} {R^{(s)}}^{\mathsf T} W^{(s)} R^{(s)} \left(\left(I - N A\right){x_{s}} + N b\right)\\
& \quad + M \left(b - A \sum_{r=1}^{p} {R^{(r)}}^{\mathsf T} W^{(r)} R^{(r)}\left(\left(I - N A\right){x_{r}} + N b\right)\right).
\end{aligned}
\]
Recalling, then, that
$
\sum_{s=1}^{p} {R^{(s)}}^{\mathsf T} W^{(s)} R^{(s)} = I,
$
it yields
\[
\begin{aligned}
F\left(x_{1}, \ldots, x_{p}\right) & = (I - M A) \sum_{s=1}^{p} {R^{(s)}}^{\mathsf T} W^{(s)} R^{(s)} \left(I - N A\right){x_{s}} + (M + N - M A N) b.
\end{aligned}
\]
Assuming, finally, two sets $X = \{x_{1}, \ldots, x_{p}\}$ and $Z = \{z_{1}, \ldots, z_{p}\}$ of vectors in $\mathbb{R}^{n}$, we have
\[
\begin{aligned}
\left|F\left(X\right) - F\left(Z\right)\right| & = \left|(I - M A) \sum_{s=1}^{p} {R^{(s)}}^{\mathsf T} W^{(s)} R^{(s)} \left(I - N A\right)\left({x_{s}} - z_{s}\right)\right|\\
& \le \sum_{s=1}^{p} \left|(I - M A) {R^{(s)}}^{\mathsf T} W^{(s)} R^{(s)} \left(I - N A\right)\right| \max_{s=1}^{p}\left(\left|x_{s} - z_{s}\right|\right).
\end{aligned}
\]
If, therefore, the condition \eqref{eq:cai_conv_cond} is satisfied, then $F$ is $p$-contracting with respect to $|.|$, which is a sufficient condition for convergence, according to the Baudet's nonlinear theory of asynchronous iterations with memory \cite{Baudet1978}.
\end{proof}

It has been shown in \cite{FrommerSzyld2001}, for several Schwarz-type preconditioners, that
one can verify both $I-MA \ge 0$ and $I-NA \ge 0$ if $A$ is an $\mathsf M$-matrix. This has been particularly pointed out in \cite{BenziEtAl2001} for $I-NA$ when $\widetilde P = \widetilde R^{\mathsf T}$ and $\widetilde A = \widetilde R A \widetilde R^{\mathsf T}$. Furthermore, \cite{FrommerSzyld2001} shows that under such conditions, one satisfies
\[
\|(I-MA)(I-NA)\|_{\infty}^{w} < 1.
\]
Restating Theorem \ref{theo:cai_conv} in the following slightly more restrictive, but more readable form can therefore be of some interest.
\begin{lemma}
\label{lem:cai_conv_mmatrix}
A two-level asynchronous iterative method \eqref{eq:cai} is convergent if
\[
\rho\left(|I-MA||I-NA|\right) < 1, \qquad N := \widetilde P \widetilde A^{-1} \widetilde R.
\]
\end{lemma}
\begin{proof}
This directly follows from the proof of Theorem \ref{theo:cai_conv}, noting that
\[
\sum_{s=1}^{p} \left|(I - M A) {R^{(s)}}^{\mathsf T} W^{(s)} R^{(s)} \left(I - N A\right)\right| \le |I-MA||I-NA|.
\]
\end{proof}
A convergence result for $\mathsf M$-matrices can now be stated.
\begin{theorem}
\label{theo:cai_conv_mmatrix}
Let $\widetilde R$ be a restriction mapping having at least one row in common with each $R^{(s)}$, $s \in \{1, \ldots, p\}$, let $\widetilde P = \widetilde R^{\mathsf T}$ and $\widetilde A = \widetilde R A \widetilde P$. Then, a two-level asynchronous Schwarz-type method \eqref{eq:cai} is convergent if $A$ is an $\mathsf M$-matrix.
\end{theorem}
\begin{proof}
By \cite[Section 8]{FrommerSzyld2001}, we have
\[
\|(I-MA)(I-NA)\|_{\infty}^{w} < 1,
\]
which implies, by Perron--Frobenius' theory (see, e.g., \cite[Corollary 6.1]{BertTsit1989}), that
\[
\rho\left(|(I-MA)(I-NA)|\right) < 1.
\]
From \cite{FrommerSzyld2001}, we also have both $I-MA \ge 0$ and $I-NA \ge 0$, and, therefore,
\[
\rho\left(|I-MA||I-NA|\right) = \rho\left(|(I-MA)(I-NA)|\right) < 1,
\]
which implies convergence, according to Lemma \ref{lem:cai_conv_mmatrix}.
\end{proof}
\begin{remark}
It should be stressed that, due to \cite{AxelKolo1994}, constructing a preconditioner for a symmetric positive definite matrix can be achieved by considering its ``diagonally compensated reduced'' counterpart, which turns out to be an $\mathsf{M}$-matrix.
\end{remark}

\subsection{Damping-based generalization}

{While Theorem \ref{theo:cai_conv_mmatrix} is restricted to $\mathsf M$-matrices and Schwarz-type preconditioners, asynchronous convergence can be more generally achieved by damping the coarse-space correction}, which, however, would deteriorate the algorithm efficiency.
{Still, the necessary amount of damping may happen to be low in practice, depending on the algorithm, on implementation aspects, on the problem and on the computing environment. In such suitable cases, it would remain the more straightforward way of achieving two-level asynchronous solvers without any particular assumption.
}
The damped asynchronous multiplicative coarse-space correction is given by
\begin{equation}
\label{eq:cai_damped}
\left \{
\begin{array}{lcl}
{\chi^{(s)}}^{k} & := & {x^{(s)}}^{\delta_{s,s}(k)}\\
&& +\ \theta R^{(s)} \widetilde P \widetilde A^{-1} \left(\widetilde b - \displaystyle\sum_{r=1}^{p} \widetilde R A {R^{(r)}}^{\mathsf T} W^{(r)} {x^{(r)}}^{\delta_{s,r}(k)}\right),\\
{x^{(s)}}^{k+1} & = & \left\{
\begin{array}{ll}
f^{(s)}\left({\chi^{(1)}}^{\delta_{s,1}(k)}, \ldots, {\chi^{(p)}}^{\delta_{s,p}(k)}\right) & \forall s \in S_{k},\\
{x^{(s)}}^{k} & \forall s \notin S_{k},
\end{array}
\right.
\end{array}
\right.
\end{equation}
where $\theta > 0$ is the damping factor. Classical correction is simply given by $\theta = 1$. Our analysis framework straightforwardly provides a general convergence proof.
\begin{theorem}
\label{theo:cai_damped_conv}
With $\theta$ sufficiently small, a two-level asynchronous iterative me-thod \eqref{eq:cai_damped} is convergent if so is its one-level counterpart \eqref{eq:ai}.
\end{theorem}
\begin{proof}
We have
\[
|I - M A| |I - \theta N A| \le |I - M A| (I + \theta |N A|),
\]
which decreases to $|I - M A|$ as $\theta$ decreases to $0$. By the continuity of the spectral radius, if, therefore, $\rho(|I - M A|) < 1$ (implied by the convergence of \eqref{eq:ai}), then we can reach with $\theta$ sufficiently small,
\[
\rho(|I - M A| (I + \theta |N A|)) = \rho(|I - M A|) + \varepsilon < 1, \qquad \varepsilon \in \mathbb{R},
\]
which implies convergence for \eqref{eq:cai_damped}, following the same development as in the proofs of Theorem \ref{theo:cai_conv} and Lemma \ref{lem:cai_conv_mmatrix}.
\end{proof}

\subsection{Centralized/distributed coarse problem}

We have seen from the explicit asynchronous model \eqref{eq:cai} that an iteration matrix $(I-MA)(I-NA)$ could not be deduced as done for classical (synchronous) two-level methods. More generally, while convergence conditions are the same for both classical parallel methods and their underlying sequential single-process version,
asynchronous convergence conditions, on the other side, are tightly related to the precise parallel scheme under consideration. In the model \eqref{eq:cai_damped}, a coarse solution is independently computed by each process. Let us now consider the general, damping-based, asynchronous formulation of the centralized/distributed scheme \eqref{eq:csi2}, which is given by
\begin{equation}
\label{eq:cai2_damped}
\left \{
\begin{array}{lcl}
\widetilde y^{k} & := & \widetilde A^{-1} \left(\widetilde b - \displaystyle\sum_{r=1}^{p} \widetilde R A {R^{(r)}}^{\mathsf T} W^{(r)} {x^{(r)}}^{\delta_{0,r}(k)}\right),\\
{\chi^{(s)}}^{k} & := & {x^{(s)}}^{\delta_{s,s}(k)} + \theta R^{(s)} \widetilde P \widetilde y^{\delta_{s,0}(k)},\\
{x^{(s)}}^{k+1} & = & \left\{
\begin{array}{ll}
f^{(s)}\left({\chi^{(1)}}^{\delta_{s,1}(k)}, \ldots, {\chi^{(p)}}^{\delta_{s,p}(k)}\right) & \forall s \in S_{k},\\
{x^{(s)}}^{k} & \forall s \notin S_{k},
\end{array}
\right.
\end{array}
\right.
\end{equation}
where a process (possibly a group of processes), say $0$, is dedicated to the evaluation (possibly parallel and asynchronous) of the coarse solution. Delay $k - \delta_{0,r}(k)$ is therefore observed on data transmission from each process $r \in \{1, \ldots, p\}$ toward the process $0$, as well as delay $k - \delta_{s,0}(k)$ on transmitting the coarse solution $\widetilde y$ from the process $0$ toward each process $s \in \{1, \ldots, p\}$. The main difference with \eqref{eq:cai} and \eqref{eq:cai_damped} is that the local sub-vector $x^{(s)}$ being corrected is no more necessarily the one involved in the right-hand side of the coarse problem. Here then, from each ${\chi^{(r)}}^{\delta_{s,r}(k)}$, $r \in \{1, \ldots, p\}$, both ${x^{(r)}}^{\delta_{r,r}(\delta_{s,r}(k))}$ and ${x^{(r)}}^{\delta_{0,r}(\delta_{r,0}(\delta_{s,r}(k)))}$ get involved, which are two different versions of the sub-vector $x^{(r)}$ used to compute ${x^{(s)}}^{k+1}$, $s \in S_{k}$. The induced iteration mapping $F^{(s)}$ thus operates on $2 p$ versions of each sub-vector $x^{(r)}$, $r \in \{1, \ldots, p\}$.
The convergence results for the shared/replicated scheme
therefore do not apply
since they were based on an iteration mapping with $p$ entries. Nevertheless, following the same approach as in the proof of Theorem \ref{theo:cai_conv}, Theorem \ref{theo:cai_damped_conv} is still achievable for the centralized/distributed scheme as well.
\begin{theorem}
\label{theo:cai2_damped_conv}
With $\theta$ sufficiently small, a two-level asynchronous iterative me-thod \eqref{eq:cai2_damped} is convergent if so is its one-level counterpart \eqref{eq:ai}.
\end{theorem}
\begin{proof}
The model \eqref{eq:cai2_damped} develops, for $s \in S_{k}$, as
\[
\begin{aligned}
{x^{(s)}}^{k+1} & = {\chi^{(s)}}^{\delta_{s,s}(k)} + M^{(s)} \left(b^{(s)} - \sum_{r=1}^{p} R^{(s)} A {R^{(r)}}^{\mathsf T} W^{(r)} {\chi^{(r)}}^{\delta_{s,r}(k)}\right),
\end{aligned}
\]
\[
\begin{aligned}
{x^{(s)}}^{k+1} & = {x^{(s)}}^{\delta_{s,s}(\delta_{s,s}(k))} + \theta R^{(s)} \widetilde P \widetilde A^{-1} \left(\widetilde b - \displaystyle\sum_{q=1}^{p} \widetilde R A {R^{(q)}}^{\mathsf T} W^{(q)} {x^{(q)}}^{\delta_{0,q}(\delta_{s,0}(\delta_{s,s}(k)))}\right)\\
& \quad + M^{(s)} \left(b^{(s)} - \sum_{r=1}^{p} R^{(s)} A {R^{(r)}}^{\mathsf T} W^{(r)}\left( {x^{(r)}}^{\delta_{r,r}(\delta_{s,r}(k))} \right.\right.\\
& \left.\left. \qquad\qquad\quad +\ \theta R^{(r)} \widetilde P \widetilde A^{-1} \left(\widetilde b - \displaystyle\sum_{q=1}^{p} \widetilde R A {R^{(q)}}^{\mathsf T} W^{(q)} {x^{(q)}}^{\delta_{0,q}(\delta_{r,0}(\delta_{s,r}(k)))}\right)\right)\right).
\end{aligned}
\]
Taking, then, $2 p$ vectors $x_{1}$ to $x_{p}$ and $x_{0,1}$ to $x_{0,p}$ in $\mathbb{R}^{n}$, the corresponding iteration mapping $F^{(s)}$ is given by
{
\[
\begin{aligned}
F^{(s)}\left(x_{1}, \ldots, x_{p},\right.\quad&\\\left. x_{0,1}, \ldots, x_{0,p}\right) & := {x_{s}}^{(s)} + \theta R^{(s)} \widetilde P \widetilde A^{-1} \left(\widetilde b - \displaystyle\sum_{q=1}^{p} \widetilde R A {R^{(q)}}^{\mathsf T} W^{(q)} {x_{0,s}}^{(q)}\right)\\
& \quad + M^{(s)} \left(b^{(s)} - \sum_{r=1}^{p} R^{(s)} A {R^{(r)}}^{\mathsf T} W^{(r)}\left( {x_{r}}^{(r)} \right.\right.\\
& \left.\left.\qquad\qquad\quad +\ \theta R^{(r)} \widetilde P \widetilde A^{-1} \left(\widetilde b - \displaystyle\sum_{q=1}^{p} \widetilde R A {R^{(q)}}^{\mathsf T} W^{(q)} {x_{0,r}}^{(q)}\right)\right)\right)\\
& = R^{(s)} {x_{s}} + \theta R^{(s)} \widetilde P \widetilde A^{-1} \widetilde R \left(b - A {x_{0,s}}\right)\\
& \quad + M^{(s)} R^{(s)} \left(b - A \sum_{r=1}^{p} {R^{(r)}}^{\mathsf T} W^{(r)}\left(R^{(r)} {x_{r}} \right.\right.\\
& \left.\left. \qquad\qquad\qquad\quad +\ \theta R^{(r)} \widetilde P \widetilde A^{-1} \widetilde R \left(b - A {x_{0,r}}\right)\right)\right).
\end{aligned}
\]
}
The induced global iteration mapping follows as
\[
\begin{aligned}
F\left(x_{1}, \ldots, x_{p}, x_{0,1}, \ldots, x_{0,p}\right) & = \sum_{s=1}^{p} {R^{(s)}}^{\mathsf T} W^{(s)} F^{(s)}\left(x_{1}, \ldots, x_{p}, x_{0,1}, \ldots, x_{0,p}\right)\\
& = \sum_{s=1}^{p} (I - M A) {R^{(s)}}^{\mathsf T} W^{(s)} R^{(s)} \left(x_{s} - \theta N A x_{0,s}\right)\\
& \quad\ +\ (M + \theta N - \theta M A N) b
\end{aligned}
\]
with $N := \widetilde P \widetilde A^{-1} \widetilde R$. Assuming, then, two sets $X = \{x_{1}, \ldots, x_{p}, x_{0,1}, \ldots, x_{0,p}\}$ and $Z = \{z_{1}, \ldots, z_{p}, z_{0,1}, \ldots, z_{0,p}\}$ of vectors in $\mathbb{R}^{n}$, we have
\[
\begin{aligned}
\left|F\left(X\right) - F\left(Z\right)\right| & = \left|\sum_{s=1}^{p} (I - M A) {R^{(s)}}^{\mathsf T} W^{(s)} R^{(s)} \left({x_{s}} - z_{s} - \theta N A\left({x_{0,s}} - z_{0,s}\right)\right)\right|\\
& \le \sum_{s=1}^{p} \left|(I - M A) {R^{(s)}}^{\mathsf T} W^{(s)} R^{(s)}\right| \left|{x_{s}} - z_{s}\right|\\
& \quad\ +\ \theta \sum_{s=1}^{p} \left|(I - M A) {R^{(s)}}^{\mathsf T} W^{(s)} R^{(s)}\right|\left|N A\right|\left|{x_{0,s}} - z_{0,s}\right|\\
& \le |I - M A| (I + \theta |N A|) \max_{s=1}^{p}\left(\left|x_{s} - z_{s}\right|, \left|x_{0,s} - z_{0,s}\right|\right).
\end{aligned}
\]
As in the proof of Theorem \ref{theo:cai_damped_conv}, with $\rho(|I - M A|) < 1$ and $\theta$ sufficiently small, one can have
\[
\rho(|I - M A| (I + \theta |N A|)) < 1,
\]
which implies convergence, as in the proof of Theorem \ref{theo:cai_conv}, due to the Baudet's contraction-based nonlinear theory \cite{Baudet1978} about asynchronous iterations with memory.
\end{proof}
\begin{remark}
The asynchronous version of the shared/replicated scheme leads to multiple different coarse right-hand-side (RHS) vectors.
From a practical point of view, applying, then, multiple different coarse solutions is very likely to introduce more global inconsistency than centralizing/distributing the coarse problem. While, therefore, Theorem \ref{theo:cai2_damped_conv} does not specifically ensure convergence for $\theta = 1$, this can reasonably be expected in practice, based on Lemma \ref{lem:cai_conv_mmatrix} and Theorem \ref{theo:cai_conv_mmatrix}, as confirmed by numerical experiments.
\end{remark}

\section{Implementation aspects}
\label{sec:impl}

\subsection{ Practical iterative model}

From the iterative models \eqref{eq:cai_damped} and \eqref{eq:cai2_damped}, the RHS vector of the coarse problem is given by
\[
\widetilde \tau :=  \widetilde R \tau, \qquad \tau := b - A \displaystyle\sum_{r=1}^{p} {R^{(r)}}^{\mathsf T} W^{(r)} {\bar x^{(r)}}
\]
with, respectively, either ${\bar x^{(r)}} := {x^{(r)}}^{\delta_{s,r}(k)}$ or ${\bar x^{(r)}} := {x^{(r)}}^{\delta_{0,r}(k)}$.
In practice, the residual vector $\tau$ is usually readily accessible in the form
\[
\tau = \sum_{s=1}^{p} {R^{(s)}}^{\mathsf T} W^{(s)} \tau^{(s)}, \qquad \tau^{(s)} := b^{(s)} - \sum_{r=1}^{p} R^{(s)} A {R^{(r)}}^{\mathsf T} W^{(r)} {\bar x^{(r)}},
\]
which therefore implies the assembly scheme
\[
\widetilde \tau = \sum_{s=1}^{p} \widetilde \tau^{(s)}, \qquad \widetilde \tau^{(s)} := {\widetilde R^{(s)^{\mathsf T}}} W^{(s)} \tau^{(s)}, \qquad {\widetilde R^{(s)^{\mathsf T}}} := \widetilde R {R^{(s)}}^{\mathsf T},
\]
where each process $s \in \{1, \ldots, p\}$ transmits its component $\widetilde \tau^{(s)}$ either to all other processes (shared/replicated approach) or to a dedicated one (centralized/distributed approach). The model \eqref{eq:cai_damped}, for instance, is therefore actually altered into
\[
\left \{
\begin{array}{lcl}
\widetilde \tau^{(s)^{k}} & := & {\widetilde R^{(s)^{\mathsf T}}} W^{(s)} \left( b^{(s)} - \displaystyle\sum_{r=1}^{p} R^{(s)} A {R^{(r)}}^{\mathsf T} W^{(r)} {x^{(r)}}^{\delta_{s,r}(k)}\right),\\
{\chi^{(s)}}^{k} & := & {x^{(s)}}^{\delta_{s,s}(k)} + \theta R^{(s)} \widetilde P \widetilde A^{-1} \displaystyle\sum_{r=1}^{p} \widetilde \tau^{(r)^{\delta_{s,r}(k)}},\\
{x^{(s)}}^{k+1} & = & \left\{
\begin{array}{ll}
f^{(s)}\left({\chi^{(1)}}^{\delta_{s,1}(k)}, \ldots, {\chi^{(p)}}^{\delta_{s,p}(k)}\right) & \forall s \in S_{k},\\
{x^{(s)}}^{k} & \forall s \notin S_{k},
\end{array}
\right.
\end{array}
\right.
\]
which leads to
\[
\begin{aligned}
{\chi^{(s)}}^{k}  & = {x^{(s)}}^{\delta_{s,s}(k)}\\&\quad +\ \theta R^{(s)} \widetilde P \widetilde A^{-1} \left(\widetilde b - \displaystyle\sum_{r=1}^{p} \sum_{q=1}^{p} \widetilde R {R^{(r)}}^{\mathsf T} W^{(r)} R^{(r)} A {R^{(q)}}^{\mathsf T} W^{(q)} {x^{(q)}}^{\delta_{r,q}(\delta_{s,r}(k))}\right)
\end{aligned}
\]
instead of
\[
{\chi^{(s)}}^{k} = {x^{(s)}}^{\delta_{s,s}(k)} + \theta R^{(s)} \widetilde P \widetilde A^{-1} \left(\widetilde b - \displaystyle\sum_{r=1}^{p} \widetilde R A {R^{(r)}}^{\mathsf T} W^{(r)} {x^{(r)}}^{\delta_{s,r}(k)}\right).
\]
We thus note that because of
{ ${x^{(q)}}^{\delta_{r,q}(\delta_{s,r}(k))}$},
which interlinks $r$ and $q$, the former expression cannot be simplified by considering, outside the sum on $q$,
\[
\widetilde R \displaystyle\sum_{r=1}^{p} {R^{(r)}}^{\mathsf T} W^{(r)} R^{(r)} A = \widetilde R A,
\]
{
which would have been possible only by considering delays of the form $\delta_{q}(\delta_{s}(k))$, resulting from
\[
\left \{
\begin{array}{lcl}
\widetilde \tau^{(s)^{k}} & = & {\widetilde R^{(s)^{\mathsf T}}} W^{(s)} \left( b^{(s)} - \displaystyle\sum_{r=1}^{p} R^{(s)} A {R^{(r)}}^{\mathsf T} W^{(r)} {x^{(r)}}^{\delta_{r}(k)}\right),\\
{\chi^{(s)}}^{k} & = & {x^{(s)}}^{\delta_{s,s}(k)} + \theta R^{(s)} \widetilde P \widetilde A^{-1} \displaystyle\sum_{r=1}^{p} \widetilde \tau^{(r)^{\delta_{s}(k)}},
\end{array}
\right.
\]
and leading to
\[
{\chi^{(s)}}^{k} = {x^{(s)}}^{\delta_{s,s}(k)} + \theta R^{(s)} \widetilde P \widetilde A^{-1} \left(\widetilde b - \displaystyle\sum_{r=1}^{p} \widetilde R A {R^{(r)}}^{\mathsf T} W^{(r)} {x^{(r)}}^{\delta_{r}(\delta_{s}(k))}\right),
\]
which exhibits a particular form of the general delay function, $\delta_{s,r}(k) = \delta_{r}(\delta_{s}(k))$.

Two practical constraints thus arise for implementing asynchronous iterations with multiplicative coarse-space correction. First, all of the processes $s \in \{1, \ldots, p\}$ have to compute $\widetilde \tau^{(s)}$ based on a same unique set of local vectors $x^{(r)}$, $r \in \{1, \ldots, p\}$, since the delay $\delta_{r}(k)$ does not depend on $s$. Second, any coarse RHS vector on a process $s$ must be constructed from all of the $\delta_{s}(k)$-th such computed components $\widetilde \tau^{(r)}$. The first mechanism corresponds to computing components $\tau^{(s)}$ of a unique global residual vector (restricted then as $\widetilde \tau^{(s)} := {\widetilde R^{(s)^{\mathsf T}}} W^{(s)} \tau^{(s)}$), while the second one consists of receiving all of the components $\widetilde \tau^{(r)}$ before solving the coarse problem. As we shall see, such synchronization mechanisms can be achieved in a non-blocking way, which preserves the asynchronous nature of the underlying iterative process. Note that the second mechanism was applied in \cite{WPouChow2019, GlusaEtAl2020} for the additive coarse-space correction case.
}

\subsection{Coarse RHS vector construction}

{

Computing a global residual vector during asynchronous iterations goes back to \cite{SavBert1996} for handling the convergence detection issue through what the authors called ``supervised termination''. Similarities with distributed snapshot~\cite{ChanLamp1985} was established in \cite{MagGBen2018b}, and the deduced simpler approach of non-blocking synchronization as a general collective operation is due to \cite{MagGBen2018c}. In the sequel, we shall denote non-blocking synchronization by ``ISYNC'', ``I'' standing for ``immediate'', similarly to the Message Passing Interface (MPI) standard.

To build a unique set of local vectors $x^{(r)}$, $r \in \{1, \ldots, p\}$, the ISYNC idea is to make each process $s \in \{1, \ldots, p\}$ stores its value $x^{(s)}$ (local snapshot) and provides its neighbors with interface data related to that particular $x^{(s)}$. Once a process $s$ has received all of its snapshot-tagged dependencies, which completes the ISYNC operation, it computes $\tau^{(s)}$ using such local and received snapshot data, which results in its local part of a unique global residual vector (see, e.g., \cite{GBenMag2020} for a more exhaustive discussion). A second ISYNC phase is then performed with $\widetilde \tau^{(s)}$ to build the coarse RHS vector.

}

Using semantic of classical two-sided communication routines of the MPI standard,
{
Algorithm \ref{algo:isync} gives an example of what could be an ISYNC collective routine.
\begin{algorithm}[htbp]
\caption{ISynchronize($x^{(s)}$, $\bar x$, dests, sources)}
\label{algo:isync}
{
\begin{algorithmic}[1]
\STATE $\bar x^{(s)}$ := $x^{(s)}$
\FOR{$r$ \textbf{in} dests}
	\STATE{Free(ISend($\bar x^{(s)}$, $r$))}
\ENDFOR
\FOR{$r$ \textbf{in} sources}
	\STATE{requests[$r$] := IRecv($\bar x^{(r)}$, $r$)}
\ENDFOR
\RETURN{requests}
\end{algorithmic}
}
\end{algorithm}
The buffer $\bar x$ contains the copied local component $x^{(s)}$ and the snapshot-tagged interface data $\bar x^{(r)}$, $r \in \{1, \ldots, p\}$, $r \ne s$. To simplify the description, interface data are taken as whole local components of the other processes $r$. Similarly to classical MPI non-blocking routines, completion can then be checked on the returned request objects. The second ISYNC phase, which does not require to save a copy of $\widetilde \tau^{(s)}$, would be given by \emph{ISynchronize($\widetilde \tau^{(s)}$, $\widetilde \tau^{}$, dests, sources)}. Copy could therefore be avoided based on whether a special parameter value is passed instead of $\widetilde \tau^{(s)}$, similarly to the \emph{MPI\_IN\_PLACE} keyword classically used in MPI collective routines. Note, furthermore, that such an ISYNC operation is indeed simply equivalent to an \emph{IAllGather(IN\_PLACE, $\widetilde \tau^{}$)}. For simplicity, the buffer $\widetilde \tau^{}$, here, is meant as containing each component $\widetilde \tau^{(r)}$, $r \in \{1, \ldots, p\}$, instead of denoting their sum as in the mathematical descriptions.
}

\subsection{Corrections frequency}

With the ISYNC-based coarse RHS vector construction, the coarse solution vector is updated only when a collective ISYNC operation completes, which potentially takes several asynchronous iterations. The same coarse solution may therefore be used several times to correct the changing fine solution. In \cite{WPouChow2019, GlusaEtAl2020} for the additive case, authors suggested to apply each coarse solution only once.
From our general models \eqref{eq:cai_damped} and \eqref{eq:cai2_damped}, no control is needed on the corrections frequency, however, applying the same correction several times may indeed raise performance concerns. The approach in \cite{WPouChow2019, GlusaEtAl2020} can be slightly generalized by introducing a parameter, $\zeta$, that defines the maximum allowed number of identical corrections. We shall indeed see from experimental results that $\zeta = 1$ does not necessarily correspond to the best choice.

\subsection{Overall implementation}

Algorithm \ref{algo:carhs}, \ref{algo:carhs2_root} and \ref{algo:carhs2_nonroot} provide implementation examples of replicated and centralized coarse solutions, { according to the previously described two-phase ISYNC approach.}
\begin{algorithm}[htbp]
\caption{ISYNC-based replicated coarse solution}
\label{algo:carhs}
{
\begin{algorithmic}[1]
\IF{phase = 0}
	\STATE reqs := ISynchronize($x^{(s)}$, $\bar x$, $\{1, \ldots, p\} \setminus \{s\}$, $\{1, \ldots, p\} \setminus \{s\}$)
	\STATE phase := 1
\ENDIF
\IF{phase = 1}
	\STATE end := TestAll(reqs)
	\IF{end}
		\STATE $\widetilde \tau^{(s)}$ := ${\widetilde R^{(s)^{\mathsf T}}} W^{(s)} \left( b^{(s)} - \displaystyle\sum_{r=1}^{p} R^{(s)} A {R^{(r)}}^{\mathsf T} W^{(r)} {\bar x^{(r)}}\right)$
		\STATE req := IAllGather(IN\_PLACE, $\widetilde \tau$)
		\STATE phase := 2
	\ENDIF
\ENDIF
\IF{phase = 2}
	\STATE end := Test(req)
	\IF{end}
		\STATE{$\widetilde x$ := $\widetilde A^{-1} \displaystyle\sum_{r=1}^{p} \widetilde \tau^{(r)}$}
		\STATE phase := 0
		\STATE nbidentcorr := 0
	\ENDIF
\ENDIF
\end{algorithmic}
}
\end{algorithm}
\begin{algorithm}[htbp]
\caption{ISYNC-based centralized coarse solution (root process $s_{0}$)}
\label{algo:carhs2_root}
{
\begin{algorithmic}[1]
\IF{phase = 0}
	\STATE reqs := ISynchronize($x^{(s)}$, $\bar x$, $\{1, \ldots, p\} \setminus \{s\}$, $\{1, \ldots, p\} \setminus \{s\}$)
	\STATE phase := 1
\ENDIF
\IF{phase = 1}
	\STATE end := TestAll(reqs)
	\IF{end}
		\STATE $\widetilde \tau^{(s)}$ := ${\widetilde R^{(s)^{\mathsf T}}} W^{(s)} \left( b^{(s)} - \displaystyle\sum_{r=1}^{p} R^{(s)} A {R^{(r)}}^{\mathsf T} W^{(r)} {\bar x^{(r)}}\right)$
		\STATE req := IGather($\widetilde \tau$, $s_{0}$)
		\STATE phase := 2
	\ENDIF
\ENDIF
\IF{phase = 2}
	\STATE end := Test(req)
	\IF{end}
		\STATE{$\widetilde x$ := $\widetilde A^{-1} \displaystyle\sum_{r=1}^{p} \widetilde \tau^{(r)}$}
		\STATE Free(IBcast($\widetilde x$, $s_{0}$))
		\STATE phase := 0
		\STATE nbidentcorr := 0
	\ENDIF
\ENDIF
\end{algorithmic}
}
\end{algorithm}
\begin{algorithm}[htbp]
\caption{ISYNC-based centralized coarse solution (non-root process)}
\label{algo:carhs2_nonroot}
{
\begin{algorithmic}[1]
\IF{phase = 0}
	\STATE reqs := ISynchronize($x^{(s)}$, $\bar x$, $\{1, \ldots, p\} \setminus \{s\}$, $\{1, \ldots, p\} \setminus \{s\}$)
	\STATE phase := 1
\ENDIF
\IF{phase = 1}
	\STATE end := TestAll(reqs)
	\IF{end}
		\STATE $\widetilde \tau^{(s)}$ := ${\widetilde R^{(s)^{\mathsf T}}} W^{(s)} \left( b^{(s)} - \displaystyle\sum_{r=1}^{p} R^{(s)} A {R^{(r)}}^{\mathsf T} W^{(r)} {\bar x^{(r)}}\right)$
		\STATE Free(IGather($\widetilde \tau^{(s)}$, $s_{0}$))
		\STATE req := IBcast($\widetilde x$, $s_{0}$)
		\STATE phase := 2
	\ENDIF
\ENDIF
\IF{phase = 2}
	\STATE end := Test(req)
	\IF{end}
		\STATE phase := 0
		\STATE nbidentcorr := 0
	\ENDIF
\ENDIF
\end{algorithmic}
}
\end{algorithm}
{
A variable \emph{nbidentcorr} (number of identical corrections) is set to 0 at the end of each new coarse resolution, and will be incremented at each subsequent correction of the fine solution, which will happen at most $\zeta$ times using that same coarse solution.
}
An overall implementation example of two-level asynchronous methods with multiplicative coarse-space correction is given by Algorithm \ref{algo:caicai2_damped}.
\begin{algorithm}[htbp]
\caption{ Two-level asynchronous iterative solver}
\label{algo:caicai2_damped}
{
\begin{algorithmic}[1]
\STATE{FineInit($x$, fine\_rcvreq, fine\_sndreq)}
\STATE{$\tau^{(s)}$ := $b^{(s)} - \displaystyle\sum_{r=1}^{p} R^{(s)} A {R^{(r)}}^{\mathsf T} W^{(r)} x^{(r)}$}
\STATE{rdcreq := IAllReduce(${\tau^{(s)}}^{\mathsf T}W^{(s)}\tau^{(s)}$, $\tau^{\mathsf T} \tau$, SUM)}
\STATE{Wait(rdcreq)}
\STATE{$\|\tau\|$ := $\displaystyle\sqrt{\tau^{\mathsf T} \tau}$}
\STATE{k := 0}
\STATE {nbidentcorr} := 0
\STATE{phase := 0}
\WHILE{$\|\tau\| > \varepsilon$ \AND k $<$ k\_max}
	\STATE{\textbf{include one of} Algorithm \ref{algo:carhs}, \ref{algo:carhs2_root}, \ref{algo:carhs2_nonroot} \COMMENT{coarse solution procedure}}
	\IF{{nbidentcorr} $< \zeta$}
		\STATE{$x^{(s)}$ := $x^{(s)} + \theta {\widetilde R^{(s)} \widetilde x}$}
		\STATE AsyncSendRecv($x$, fine\_rcvreq, fine\_sndreq)
		\STATE {nbidentcorr} := {nbidentcorr} + 1
	\ENDIF
	\STATE{$x^{(s)}$ := $f^{(s)}\left({x^{(1)}}, \ldots, {x^{(p)}}\right)$}
	\STATE AsyncSendRecv($x$, fine\_rcvreq, fine\_sndreq)
	\STATE{$\tau^{(s)}$ := $b^{(s)} - \displaystyle\sum_{r=1}^{p} R^{(s)} A {R^{(r)}}^{\mathsf T} W^{(r)} x^{(r)}$}
	\IF{Test(rdcreq)}
		\STATE{$\|\tau\|$ := $\displaystyle\sqrt{\tau^{\mathsf T} \tau}$}
		\STATE{rdcreq := IAllReduce(${\tau^{(s)}}^{\mathsf T}W^{(s)}\tau^{(s)}$, $\tau^{\mathsf T} \tau$, SUM)}
		\STATE{k := k + 1}
	\ENDIF
\ENDWHILE
\end{algorithmic}
}
\end{algorithm}
Communication requests for the underlying one-level algorithm are initialized following Algorithm \ref{algo:ainit}, which is based on ideas from \cite{MagGBen2018c}.
\begin{algorithm}[htbp]
\caption{FineInit($x$, rcvreq, sndreq)}
\label{algo:ainit}
{
\begin{algorithmic}[1]
\FOR{$r \in \{1, \ldots, s-1, s+1, \ldots, p\}$}
	\FOR{$j \in \{1, \ldots, \text{nrcvreqs\_per\_neighb}\}$}
		\STATE{rcvreq[$r$][$j$] := IRecv($x^{(r)}$, $r$)}
	\ENDFOR
	\STATE{sndreq[$r$] := REQUEST\_NULL}
\ENDFOR
\end{algorithmic}
}
\end{algorithm}
To take the most out of asynchronous iterations, one has to allow for updating interface dependencies during computation phases. This can be achieved in the two-sided communication pattern by keeping several reception requests active on the same buffer $x^{(r)}$. The maximum number of simultaneously active requests per neighbor process (variable \emph{nrcvreqs\_per\_neighb}) is ideally set according to both the network speed and the compute nodes heterogeneity, or can even be finely tuned for each neighbor according to a speed ratio associated to each communication link (in case such information can be inferred). Dependencies are then updated as described by Algorithm \ref{algo:acom}.
\begin{algorithm}[htbp]
\caption{AsyncSendRecv($x$, rcvreq, sndreq)}
\label{algo:acom}
{
\begin{algorithmic}[1]
\FOR{$r \in \{1, \ldots, s-1, s+1, \ldots, p\}$}
	\IF{Test(sndreq[$r$])}
		\STATE{sndreq[$r$] := ISSend($x^{(s)}$, $r$)}
	\ENDIF
\ENDFOR
\FOR{$r \in \{1, \ldots, s-1, s+1, \ldots, p\}$}
	\FOR{$j \in \{1, \ldots, \text{nrcvreqs\_per\_neighb}\}$}
		\IF{Test(rcvreq[$r$][$j$])}
			\STATE{rcvreq[$r$][$j$] := IRecv($x^{(r)}$, $r$)}
		\ENDIF
	\ENDFOR
\ENDFOR
\end{algorithmic}
}
\end{algorithm}
For each request, completion is checked, upon which a new one is triggered. Non-blocking synchronous message sending (\emph{ISSend}) avoids internal MPI buffering, therefore, it allows for updating the outgoing buffer as soon as a new $x^{(s)}$ is computed, even while a sending request initially triggered with a previous value is still active. This compensates, to some extent, the need for multiple simultaneously active reception requests. We refer to \cite{MagGBen2018c} for more insight about asynchronous programming options, along with experimental investigation. Compared to \cite{MagGBen2018c} though, note here the simpler approach for reliably use the classical loop stopping criterion, { $\|\tau\| > \varepsilon$}, which is due to the more recent work \cite{GBenMag2020}.

\section{Experimental results}
\label{sec:exp}

\subsection{Preliminaries}

Similarly to \cite{GlusaEtAl2020}, we experimentally addressed the Poisson's equation,
\[
- \Delta u = g,
\]
on the three-dimensional domain $[0, 1]^{3}$, with a uniform source { $g = 4590$}, and with $u = 0$ on the boundary.
{ The subdomains were defined by partitioning the interval $[0, 1]$ on each of the three dimensions (see Table \ref{tab:exp:cfg}), then,} the discrete corresponding problem was obtained in each subdomain by P1 finite-element approximation.
The size of the local problems was kept constant at approximately 20,000 unknowns while varying the number of subdomains (weak scaling).
{
The discrete overlap was also kept constant at 2 mesh steps (see Table \ref{tab:exp:cfg}).
}
The coarse problem was algebraically derived as $\widetilde A = \widetilde R A \widetilde R^{\mathsf T}$, according to Theorem \ref{theo:cai_conv_mmatrix}, with one unknown per subdomain
{
and $\widetilde R$ being of the form
\[
\widetilde R =
\begin{bmatrix}
1 & \cdots & 1 & 0 & \cdots & 0 & 0 & \cdots & 0\\
0 & \cdots & 0 & 1 & \cdots & 1 & 0 & \cdots & 0\\
0 & \cdots & 0 & 0 & \cdots & 0 & 1 & \cdots & 1
\end{bmatrix},
\]
for an example with 3 subdomains. Unless otherwise stated,
}
the centralized scheme (Algorithm \ref{algo:carhs2_root}, \ref{algo:carhs2_nonroot}) was considered with both local fine and global coarse problems being solved through LU factorization.

Each node of the compute cluster consisted of 160 GB of RAM and two 20-core CPUs at 2.1 GHz (40 cores per node). Each core was mapped with exactly one MPI process and each MPI process with one subdomain. The compute cluster is connected through an Omni-Path Architecture (OPA) network at 100 Gbit/s.

In the following, we compare overall compute times of the different solvers, final residual errors $\|\tau\| = \|b - Ax\|$ (synchronously computed after the iterations loop) and numbers of iterations (average numbers in case of asynchronous iterations).
{ The tolerance for the stopping criterion was set to $\varepsilon = 10^{-6}$. The different parallel configurations are summarized in Table \ref{tab:exp:cfg}
\begin{table}[htbp]
\caption{Experiments configuration.}
\label{tab:exp:cfg}
\centering
{\footnotesize
\begin{tabular}{ccc}
\hline\noalign{\smallskip}
Domain & Subdomains & Cluster\\
\noalign{\smallskip}\hline\noalign{\smallskip}\noalign{\smallskip}
\begin{tabular}{c}
$n$\\
\noalign{\smallskip}\hline\noalign{\smallskip}
$80^{3}$\\
$160^{3}$\\
$320^{3}$\\
\end{tabular}
&
\begin{tabular}{ccc}
$p$ & Partitioning & Overlap\\
\noalign{\smallskip}\hline\noalign{\smallskip}
25 & $5\times 5\times 1$ & 2\\
200 & $5\times 8\times 5$ & 2\\
1600 & $10\times 16\times 10$ & 2\\
\end{tabular}
&
\begin{tabular}{ccc}
CPU cores & CPUs & Nodes\\
\noalign{\smallskip}\hline\noalign{\smallskip}
25 & 2 & 1\\
200 & 10 & 5\\
1600 & 80 & 40\\
\end{tabular}\\
\noalign{\smallskip}\hline
\end{tabular}
}
\end{table}

}

\subsection{Optimal implementation}

\begin{figure}[htbp]
\centering
\includegraphics[width=0.45\textwidth]{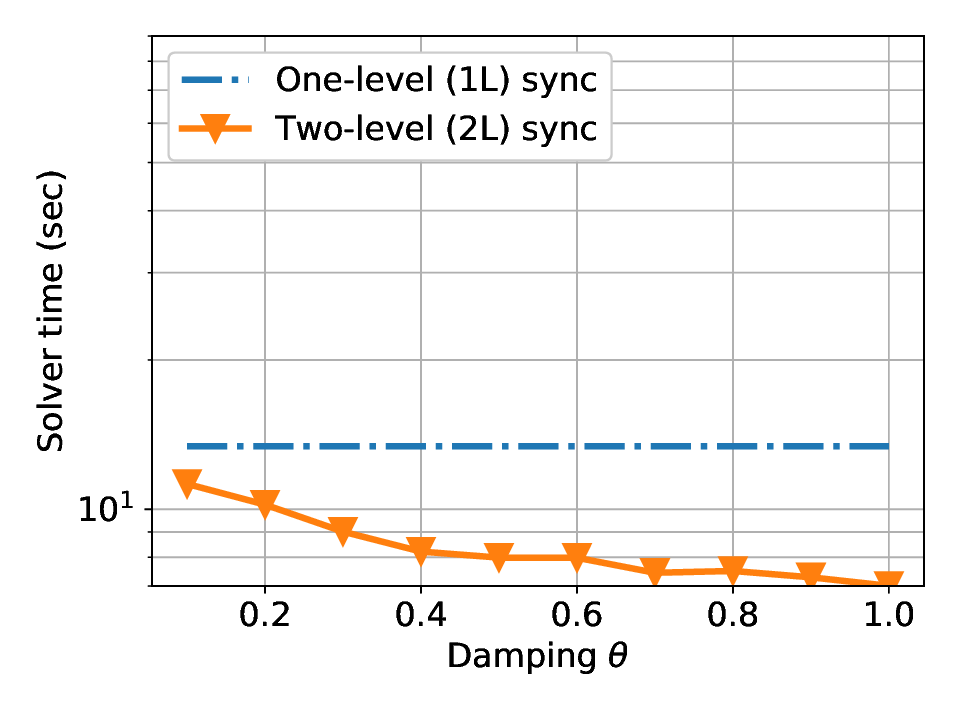}
\includegraphics[width=0.45\textwidth]{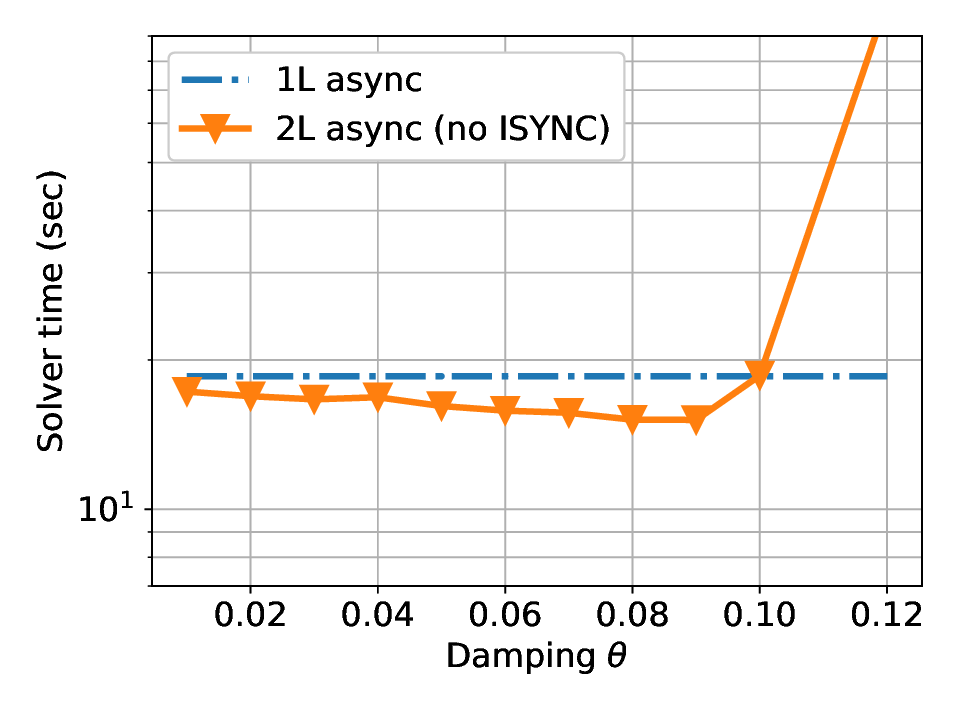}
\includegraphics[width=0.45\textwidth]{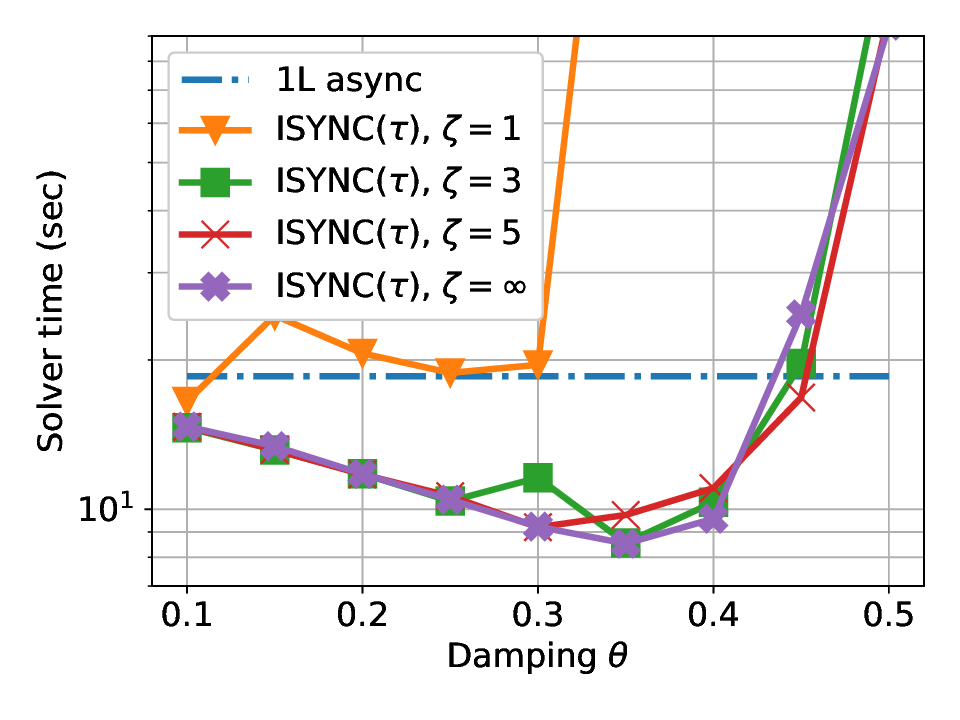}
\includegraphics[width=0.45\textwidth]{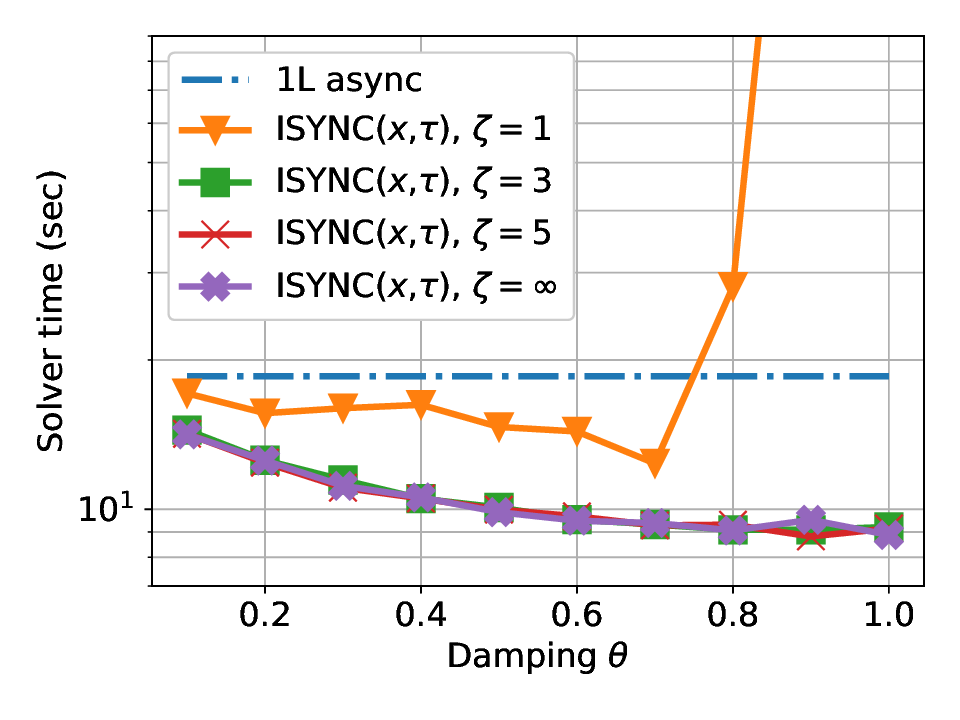}
\caption{Performance according to the correction damping, $\theta$, and the maximum allowed number of successive identical corrections, $\zeta$. Experiments with 25 CPU cores.}
\label{fig:exp:damp25}
\end{figure}

{
Algorithm \ref{algo:carhs2_root}, \ref{algo:carhs2_nonroot} describe a two-level asynchronous solver based on two non-blocking synchronization phases, successively applied to the solution and the residual vectors, $x$ and $\tau$, which we found to be the only way to match the theoretical iterative model (as a special case though). Let us refer to such an implementation approach as ISYNC($x$,$\tau$). In the additive correction form \cite{WPouChow2019, GlusaEtAl2020}, however, only ISYNC($\tau$) was considered. One might therefore wonder, for the multiplicative case, to which extent ISYNC is needed in practice. A comparative example is illustrated in Figure \ref{fig:exp:damp25}.
}
Two parameters have to be considered, the correction damping $\theta$ from the model \eqref{eq:cai2_damped} and, from Algorithm \ref{algo:caicai2_damped}, $\zeta$, the maximum number of successive identical corrections which is allowed. Strictly applying the general asynchronous model \eqref{eq:cai2_damped} corresponds to taking $\zeta = \infty$ (no bounding).

{
It clearly appears that with less of ISYNC, more damping was needed to make the solver converge. Significant differences also came out with whether $\zeta = 1$ or not, showing a particularly poor performance for $\zeta = 1$.
}
Considering
both the (average) number of iterations, $k$, and the number of computed coarse solutions, $c$,
the (average) number of successive identical corrections, $k/c$, was varying
{
between 3 and 4 for ISYNC($x$,$\tau$) and between 2 and 3 for ISYNC($\tau$), independently of $\theta$. This possibly explains the performance invariability with bounds $\zeta \in \{3, 5, \infty\}$. Note that ISYNC($x$,$\tau$) with $\zeta = \infty$ remained fully effective at 1600 processor cores where variations of $k/c$ increased to the interval $[6, 7]$. Implementation with $\zeta > 1$ or even without $\zeta$ ($\zeta = \infty$) therefore appears to be possible here, contrarily to \cite{WPouChow2019, GlusaEtAl2020}. Note, finally, on Figure \ref{fig:exp:damp25} that despite a pronounced damping $\theta = 0.35$, ISYNC($\tau$) shows competitiveness with ISYNC($x$,$\tau$). Unfortunately, it drastically deteriorated when switching to 200 and 1600 cores. The theoretical necessity for an ISYNC($x$,$\tau$)-based implementation is thus confirmed in the considered practical framework.
}

\subsection{Scaling}

\begin{figure}[htbp]
\centering
\includegraphics[width=0.45\textwidth]{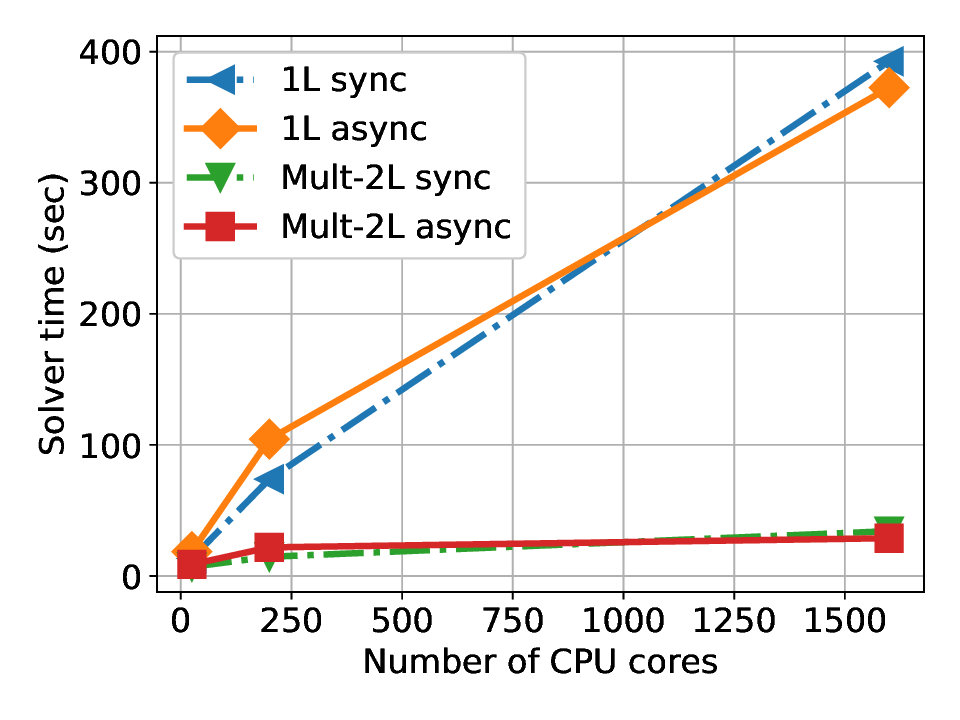}
\includegraphics[width=0.45\textwidth]{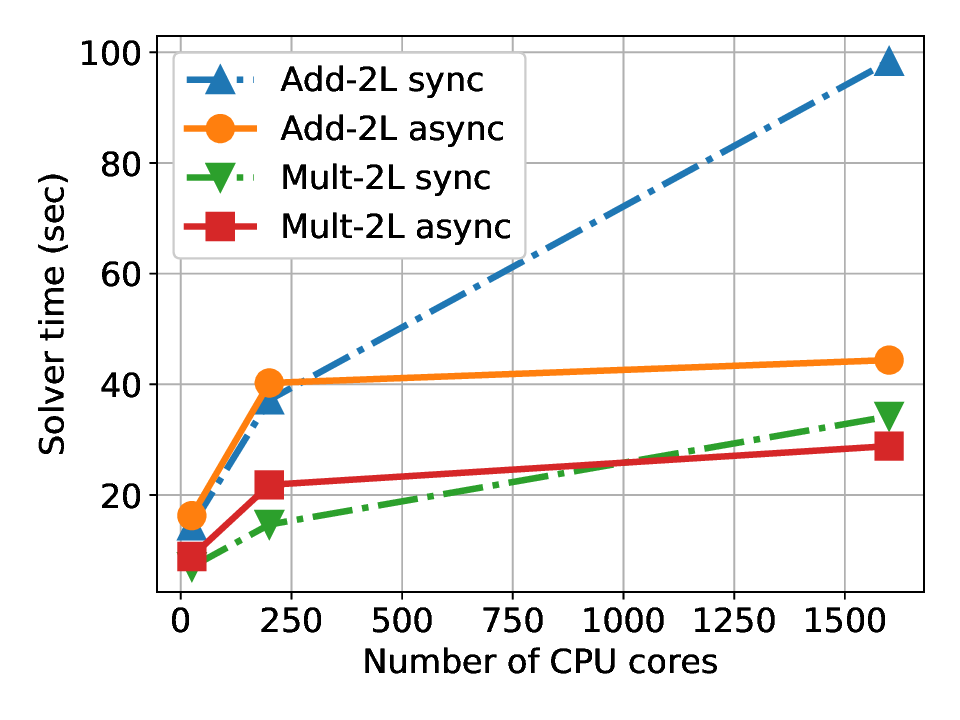}
\caption{Weak scaling performance. For the two-level asynchronous methods, the ISYNC($x$,$\tau$)-based implementation is considered with $\theta = 1$ and $\zeta = \infty$.}
\label{fig:exp:scaling}
\end{figure}

Table \ref{tab:exp:scaling} and Figure \ref{fig:exp:scaling} show the parallel computational performance of the one-level and two-level methods while varying the number of CPU cores but keeping about 20,000 fine-unknowns per core and 1 coarse-unknown per core.
\begin{table}[htbp]
\caption{Weak scaling performance. For the two-level asynchronous methods, the ISYNC($x$,$\tau$)-based implementation is considered with $\theta = 1$ and $\zeta = \infty$. Solver times (T) are in seconds.}
\label{tab:exp:scaling}
\centering
{\footnotesize
\begin{tabular}{ccc}
\hline\noalign{\smallskip}
& Synchronous & Asynchronous \\
\noalign{\smallskip}\hline\noalign{\smallskip}\noalign{\smallskip}
\begin{tabular}{cc}
$p$ & Coarse\\
\noalign{\smallskip}\hline\noalign{\smallskip}
25
&
\begin{tabular}{c}
--\\
Add\\
Mult
\end{tabular}\\
\noalign{\smallskip}\hline\noalign{\smallskip}
200
&
\begin{tabular}{c}
--\\
Add\\
Mult
\end{tabular}\\
\noalign{\smallskip}\hline\noalign{\smallskip}
1600
&
\begin{tabular}{c}
--\\
Add\\
Mult
\end{tabular}\\
\end{tabular}
&
\begin{tabular}{ccc}
T & $k$ & $\|b - Ax\|$\\
\noalign{\smallskip}\hline\noalign{\smallskip}
\begin{tabular}{c}
13\\
14\\
7
\end{tabular}
&
\begin{tabular}{c}
475\\
501\\
213
\end{tabular}
&
\begin{tabular}{c}
9.9E-07\\
9.8E-07\\
9.8E-07
\end{tabular}\\
\noalign{\smallskip}\hline\noalign{\smallskip}
\begin{tabular}{c}
74\\
37\\
15
\end{tabular}
&
\begin{tabular}{c}
1553\\
690\\
222
\end{tabular}
&
\begin{tabular}{c}
9.9E-07\\
9.9E-07\\
9.5E-07
\end{tabular}\\
\noalign{\smallskip}\hline\noalign{\smallskip}
\begin{tabular}{c}
393\\
98\\
34
\end{tabular}
&
\begin{tabular}{c}
5669\\
770\\
224
\end{tabular}
&
\begin{tabular}{c}
9.9E-07\\
9.9E-07\\
9.8E-07
\end{tabular}\\
\end{tabular}
&
\begin{tabular}{ccc}
T & $k$ & $\|b - Ax\|$\\
\noalign{\smallskip}\hline\noalign{\smallskip}
\begin{tabular}{c}
19\\
16\\
9
\end{tabular}
&
\begin{tabular}{c}
770\\
668\\
329
\end{tabular}
&
\begin{tabular}{c}
7.1E-07\\
7.0E-07\\
9.0E-07
\end{tabular}\\
\noalign{\smallskip}\hline\noalign{\smallskip}
\begin{tabular}{c}
104\\
40\\
22
\end{tabular}
&
\begin{tabular}{c}
2472\\
909\\
442
\end{tabular}
&
\begin{tabular}{c}
7.7E-07\\
8.3E-07\\
5.0E-07
\end{tabular}\\
\noalign{\smallskip}\hline\noalign{\smallskip}
\begin{tabular}{c}
373\\
44\\
29
\end{tabular}
&
\begin{tabular}{c}
8620\\
987\\
542
\end{tabular}
&
\begin{tabular}{c}
8.3E-07\\
8.0E-07\\
4.8E-07
\end{tabular}\\
\end{tabular}\\
\noalign{\smallskip}\hline
\end{tabular}
}
\end{table}
{
As we can see, in terms of computational time, the ISYNC-asynchronous multiplicative coarse-space correction provided scaling performances comparable to its classical synchronous counterpart, which it even outperformed at 1600 processor cores. Furthermore, considering its asynchronous nature, and compared to its one-level counterpart, its slightly increasing average number of iterations can be considered as being nearly constant, which strengthens the scaling ability of the approach. Its additive counterpart (based here on ISYNC($x$,$\tau$)) showed quite similar good performances and even provided a more pronounced synchronous-to-asynchronous speedup. Still, the solver with multiplicative correction remained faster as in the synchronous case.
}

{
Solving the coarse problem using a direct method can become particularly time-consuming as the number of subdomains increases. It is therefore usual to resort to a Krylov method instead. While we did not have access to such higher levels of parallelism, Table \ref{tab:exp:scaling:cg} illustrates the scaling behavior with both local and coarse solvers based on the conjugate gradient (CG) method. In particular, the replicated coarse-problem scheme was considered here.
\begin{table}[htbp]
\caption{Weak scaling performance with CG local/coarse solvers with a tolerance of $10^{-9}$. The coarse problem is independently solved on each process (replicated scheme). For the two-level asynchronous methods, the ISYNC($x$,$\tau$)-based implementation is considered with $\theta = 1$ and $\zeta = \infty$. Overall solver times (T) are in seconds.}
\label{tab:exp:scaling:cg}
\centering
{\footnotesize
\begin{tabular}{ccc}
\hline\noalign{\smallskip}
& Synchronous & Asynchronous\\
\noalign{\smallskip}\hline\noalign{\smallskip}\noalign{\smallskip}
\begin{tabular}{cc}
$p$ & Coarse\\
\noalign{\smallskip}\hline\noalign{\smallskip}
25
&
\begin{tabular}{c}
--\\
Add\\
Mult
\end{tabular}\\
\noalign{\smallskip}\hline\noalign{\smallskip}
200
&
\begin{tabular}{c}
--\\
Add\\
Mult
\end{tabular}\\
\noalign{\smallskip}\hline\noalign{\smallskip}
1600
&
\begin{tabular}{c}
--\\
Add\\
Mult
\end{tabular}\\
\end{tabular}
&
\begin{tabular}{ccc}
T & $k$ & $\|b - Ax\|$\\
\noalign{\smallskip}\hline\noalign{\smallskip}
\begin{tabular}{c}
14\\
20\\
9
\end{tabular}
&
\begin{tabular}{c}
475\\
501\\
213
\end{tabular}
&
\begin{tabular}{c}
9.9E-07\\
9.8E-07\\
9.8E-07
\end{tabular}\\
\noalign{\smallskip}\hline\noalign{\smallskip}
\begin{tabular}{c}
79\\
40\\
14
\end{tabular}
&
\begin{tabular}{c}
1553\\
690\\
222
\end{tabular}
&
\begin{tabular}{c}
9.9E-07\\
9.9E-07\\
9.5E-07
\end{tabular}\\
\noalign{\smallskip}\hline\noalign{\smallskip}
\begin{tabular}{c}
280\\
151\\
51
\end{tabular}
&
\begin{tabular}{c}
5780\\
770\\
224
\end{tabular}
&
\begin{tabular}{c}
9.9E-07\\
9.9E-07\\
9.8E-07
\end{tabular}\\
\end{tabular}
&
\begin{tabular}{ccc}
T & $k$ & $\|b - Ax\|$\\
\noalign{\smallskip}\hline\noalign{\smallskip}
\begin{tabular}{c}
19\\
23\\
12
\end{tabular}
&
\begin{tabular}{c}
744\\
687\\
341
\end{tabular}
&
\begin{tabular}{c}
8.5E-07\\
7.2E-07\\
8.1E-07
\end{tabular}\\
\noalign{\smallskip}\hline\noalign{\smallskip}
\begin{tabular}{c}
109\\
43\\
28
\end{tabular}
&
\begin{tabular}{c}
2674\\
936\\
592
\end{tabular}
&
\begin{tabular}{c}
7.3E-07\\
8.7E-07\\
7.3E-07
\end{tabular}\\
\noalign{\smallskip}\hline\noalign{\smallskip}
\begin{tabular}{c}
358\\
67\\
37
\end{tabular}
&
\begin{tabular}{c}
11376\\
1419\\
625
\end{tabular}
&
\begin{tabular}{c}
7.7E-07\\
4.1E-07\\
6.9E-07
\end{tabular}\\
\end{tabular}\\
\noalign{\smallskip}\hline
\end{tabular}
}
\end{table}
The main drawback of such a scheme (implemented as is) is that the coarse communication graph is complete (all-to-all transmissions of $\widetilde \tau^{(s)}$), which implies $\mathcal O(p^{2})$ transfers instead of $\mathcal O(p)$ in the centralized case. This was particularly impacting at 1600 CPU cores, for which one can notice a slowdown of the solvers, compared to Table \ref{tab:exp:scaling}. Nevertheless, the overall scaling behaviors were preserved when switching to a CG local/coarse solver, as shown by Figure \ref{fig:exp:scaling:cg}.
\begin{figure}[htbp]
\centering
\includegraphics[width=0.45\textwidth]{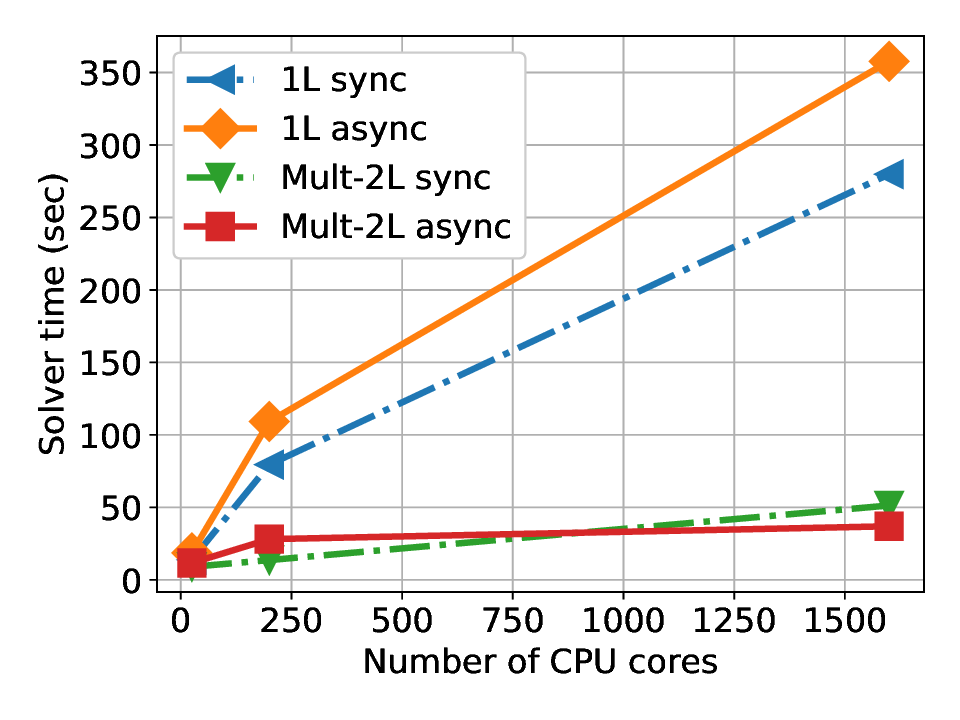}
\includegraphics[width=0.45\textwidth]{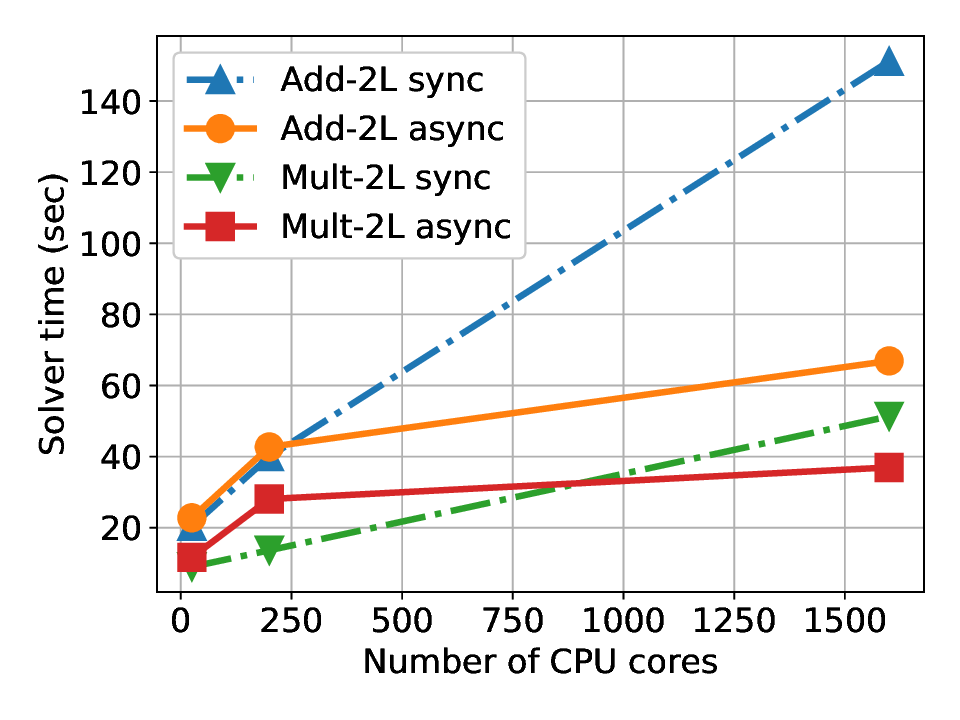}
\caption{Weak scaling performance with CG local/coarse solvers with a tolerance of $10^{-9}$. The coarse problem is independently solved on each process (replicated scheme). For the two-level asynchronous methods, the ISYNC($x$,$\tau$)-based implementation is considered with $\theta = 1$ and $\zeta = \infty$.}
\label{fig:exp:scaling:cg}
\end{figure}

\subsection{Load imbalance}

Asynchronous methods are generally expected to be less sensitive to load imbalance than their synchronous counterparts. In the two-level case, however, the coarse solver still involves synchronization (even though it is in a non-blocking way), which possibly invalidate such expectations when considering how impacted would be the benefit from coarse-space correction. What could even worsen the asynchronous solver is that an imbalanced load implies an imbalanced number of identical corrections, hence, potentially exacerbated overcorrections for underloaded processes. In such cases, a finite bound $\zeta$ could be needed back.

To simulate load imbalance in our experimental framework, both the local and the coarse solvers are slowed down by a given factor, for different subsets of processes. Note that the overall solver time includes communication time, therefore, a $2\times$-slower local solver does not imply a $2\times$-slower global synchronous solver, especially when non-blocking reception routines are used to allow for partial overlap between computation and communication. Given a maximum slowdown factor $m \in \mathbb{N}$, $m \le p$, the set of $p$ processes was partitioned into $m$ subsets $P_{1}$ to $P_{m}$, and each process $s \in \{1, \ldots, p\}$ was slowed down by a factor $i \in \{1, \ldots, m\}$ when $s \in P_{i}$. Under such simulated heterogeneous cluster, performances of the synchronous and the ISYNC-asynchronous two-level solvers can be compared in Figure \ref{fig:exp:imbal}.
\begin{figure}[htbp]
\centering
\includegraphics[width=0.45\textwidth]{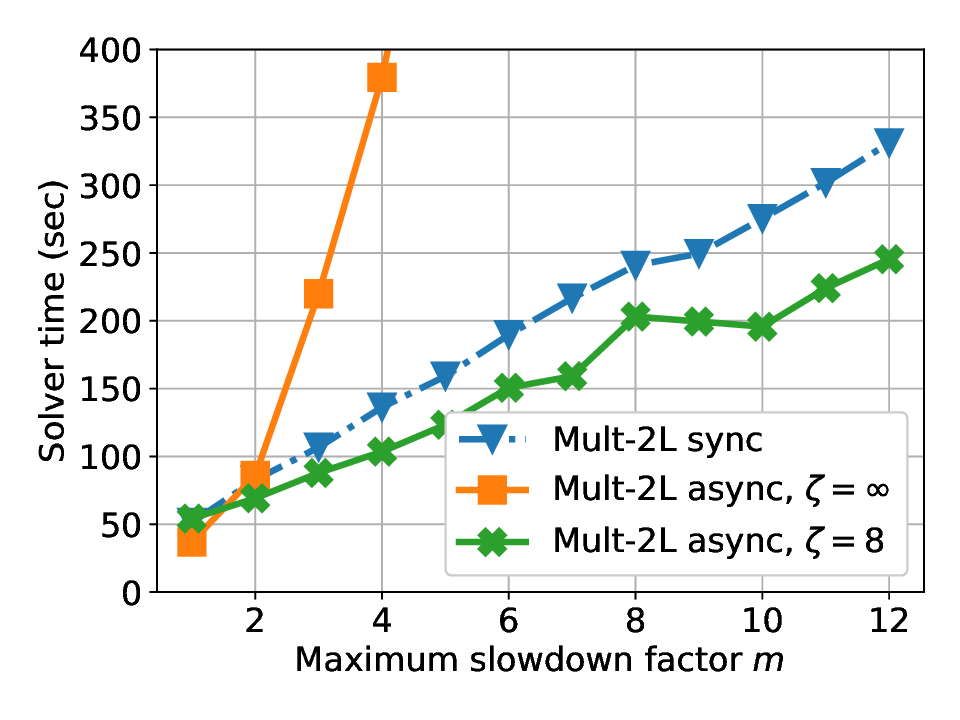}
\includegraphics[width=0.45\textwidth]{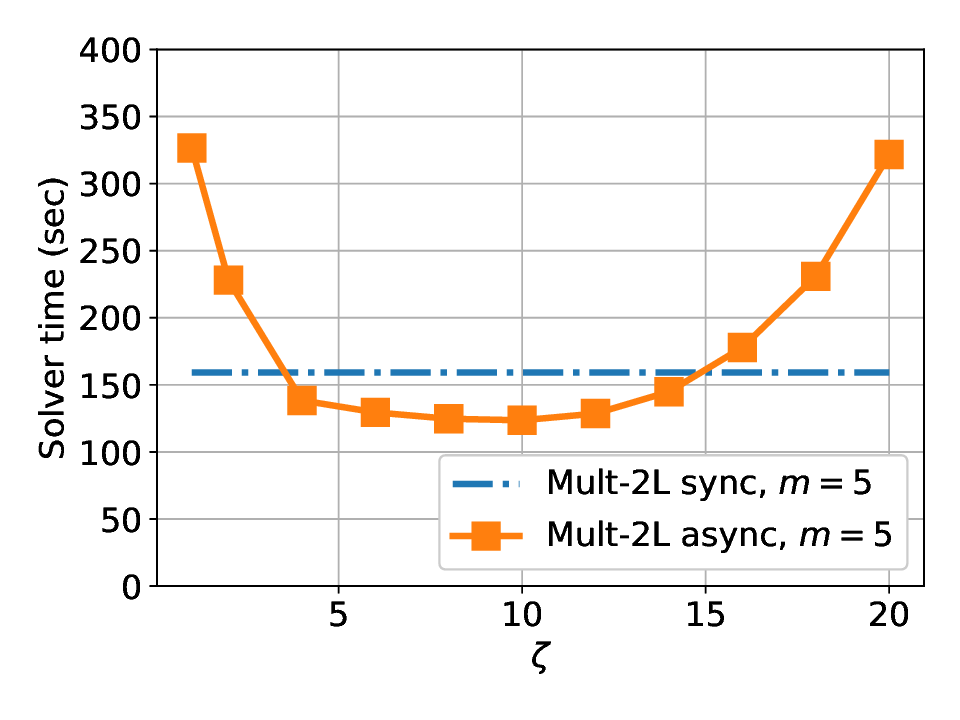}
\caption{Performance with 1600 heterogeneously slowed down processes. CG local/coarse solvers with a tolerance of $10^{-9}$ and replicated coarse-problem are considered. For the two-level asynchronous method, the ISYNC($x$,$\tau$)-based implementation is considered with $\theta = 1$.}
\label{fig:exp:imbal}
\end{figure}
As expected when increasing the heterogeneity of the cluster, the ISYNC-asynchronous solver with $\zeta = \infty$ drastically deteriorated, but quite surprisingly early. In the natural case (no simulated slower process), we observed an average number of successive identical corrections varying between $7$ and $8$. It then turned out here that taking $\zeta = 8$ allowed for retrieving better synchronous-to-asynchronous speedups. Still, such results show that it remains difficult to ensure one optimal bound for all computational environments.

}

\section{Conclusion}
\label{sec:concl}

Despite the early introduction of asynchronous iterative methods in 1969, the development and assessment of numerically efficient asynchronous preconditioners is still a big challenge. Following the recent achievement of asynchronous additive coarse-space correction, in this paper, we modeled and analyzed asynchronous multiplicative coarse-space correction.

{
In order to reach full effectiveness in practice, and even match the theoretical model, global residual vector needs to be considered instead of simply gathering available, possibly inconsistent, local components. Non-blocking synchronization techniques originally used only to solve the convergence detection issue are therefore now essential in the core implementation of asynchronous domain decomposition methods. Such combinations of asynchronous and synchronous features possibly constitute the key ingredients for achieving efficiency at very large scales.

Still, non-blocking calculation of global residual vector induces several corrections using same coarse solution, which may cause counterproductive overcorrections, especially in heterogeneous environments. Limitation of such overcorrections may therefore be needed, however, completely avoiding them does not provide the best performances. Determining an optimal bound is therefore still an issue in case of imbalanced load.
}

\section*{Acknowledgement}

The paper has been prepared with the support of the French national program LEFE/INSU and the project ADOM (M\'ethodes de d\'ecompo-sition de domaine asynchrones) of the French National Research Agency (ANR).

\bibliography{ref}
\bibliographystyle{abbrv}

\end{document}